\documentclass[12pt]{article}
\usepackage{amsmath, amsthm, amssymb}

\newtheorem{Thm}{Theorem}[section]
\newtheorem{Lem}{Lemma}[section]
\newtheorem{Prop}{Proposistion}[section]
\newtheorem{Cor}{Corollary}[section]

\begin{document}

\title{A nested embedding theorem for Hardy-Lorentz spaces with
applications to coefficient multiplier problems}

\author{
Marc Lengfield\\
{\em Department of Mathematics}\\
{\em Florida State University}\\
{\em Tallahassee, FL 32306, USA}\\
and\\
{\em Department of Mathematics}\\
{\em Western Kentucky University}\\
{\em Bowling Green, KY 42101, USA}\\
\\
Marc.Lengfield@wku.edu}

\date{\today}

\maketitle

\begin{abstract}
We prove a nested embedding theorem for Hardy-Lorentz spaces
and use it to find coefficient multiplier spaces of certain non-locally convex
Hardy-Lorentz spaces into various target spaces such as Lebesgue sequence
spaces, other Hardy spaces, and analytic mixed norm spaces.
\end{abstract}

\noindent
{\bf AMS Subject Classification:}


\section{Introduction}
In this paper we characterize coefficient multipliers between 
certain types of analytic function spaces on the open unit disk. 
We are primarily concerned with multipliers having one of the
non-locally convex Hardy-Lorentz spaces \(H^{p,q}\,,\) \(0< p < 1\,,\) 
\(0 < q \leq \infty\,,\) for the domain space. For such multipliers 
we will consider a variety of target spaces including Lebesgue
sequence spaces, other Hardy spaces, and various analytic
function spaces of mixed norm type. Our method depends upon
a nested embedding theorem for Hardy-Lorentz spaces 
(Theorem 4.1) obtained through interpolation from embedding 
theorems of Hardy and Littlewood and of Flett. Thus, the 
strategy is to trap  \(H^{p,q}\) between a pair of mixed norm
spaces of Bergman-type and then deduce multiplier results for
\(H^{p,q}\) from corresponding known multiplier results for the
endpoint spaces.
 The paper is organized as follows. In Section 2 we define the 
Hardy-Lorentz spaces and the analytic mixed norm spaces. Also
included in this section are some results from interpolation, 
fractional calculus, and \(H^p\)- theory needed for the sequel.
Our primary references for Lorentz spaces and Hardy spaces 
are [5] and [10], respectively. Section 3 covers preliminary
material on coefficient multipliers. In Section 4 we state and
prove the embedding theorem for  \(H^{p,q}\). We then indicate
how this theorem may be used to obtain the duality results of 
[33]. In addition, we determine the Abel dual of  \(H^{p,\infty}\). 
In Sections 5, 6, and 7, respectively, we find multipliers of 
\(H^{p,q}\) into the Lebesgue sequence spaces \(\ell^{s}\,,\) into 
mixed norm spaces of Bergman-type, and into certain Hardy
spaces. In Section 8 we discuss the case when the target 
space is an analytic Lipschitz or Zygmund space, a Bloch
space, or BMOA.

 Throughout the paper \(\mathbb {D}\) will denote the open unit 
disk in the complex plane and \(\mathbb {T}\) will denote its 
boundary. The symbol \(H(\mathbb {D})\) is used to denote the
space of analytic functions on \(\mathbb {D}\). If \(X\) and \(Y\) are
toplogical spaces with \(X \subset Y\) we write 
\(X\hookrightarrow Y\) to indicate continuous inclusion. All vector
spaces are assumed to be complex. By a Frechet space we 
mean a locally convex F-space. If \(E\) is a topological vector 
space then \(E^*\) denotes the topological dual space of \(E\)
consisting of all continuous linear functionals on \(E\). The 
symbol \(A\sim B\) is used to indicate the existence of absolute
positive constants \(C_j\,, j = 1,2\,,\) such that \(C_1\leq A/B \leq C_2\).
\setcounter{equation}{0}
\section{Hardy-Lorentz Spaces and Mixed Norm\\ Spaces}

Let \(m\) denote normalized Lebesgue measure on \(\mathbb {T}\) 
and let \(L^0(m)\) be the space of complex-valued Lebesgue 
measurable functions on \(\mathbb {T}\). For \(f\in L^0(m)\) and \(s \geq 0\), we  
write \(\lambda _f(s) = m(\{z\in \mathbb {T}: |f(z)|>s\}\) for the distribution function
and \(f^*(s) = \inf(\{t \geq 0 :\ \lambda_f(t)\leq s\}\) for the decreasing rearrangement of \(|f|\),
each taken with respect to \(m\). 
Let \(0 < p,q \leq \infty \). For the reader's convenience, we recall the definition of the Lorentz spaces
\(L^{p,q}(m)\). The Lorentz functional
\(||\cdot ||_{p,q}\) is defined at \(f\in L^0(m)\) by
\(|| f ||_{p,q} = (\int _{0}^{1} [f^*(s)s^{\frac{1}{p}}]^q\, \frac{ds}{s})^{1/q}\) for \(0 < q < \infty\) and 
\(|| f ||_{p,\infty }= \sup_{s\geq 0}[f^*(s)s^{\frac{1}{p}}]\).
The corresponding Lorentz space is \(L^{p,q}(m) = \{f\in L^0(m): || f ||_{p,q} < \infty \}\).
Since \(|| f ||_{p,p} = || f ||_p\,,\) where \(|| f ||_p\)
denotes the standard \(L^p\)-functional on \(L^0(m)\), the Lorentz spaces form a 2-parameter 
array $\{(L^{p,q}(m), ||\cdot ||_{p,q})\}_{0 < p, q\leq \infty }$ of quasi-Banach
spaces containing the Lebesgue space scale \(\{(L^p(m), ||f||_p)\}_{0 < p\leq \infty }\) as the main diagonal. 
Inclusions among the Lorentz spaces are given by

\begin{equation}\label{2.1}
L^{p,q}(m) \hookrightarrow L^{p,r}(m)\,,\\ 0< p\leq  \infty\,, 0< q \leq r \leq \infty\,,
\end{equation}
\noindent and, since \(m(\mathbb{T}) < \infty\),

\begin{equation}\label{2.2}
L^{r,s}(m) \hookrightarrow L^{p,q}(m)\,,\\ 0< p < r \leq \infty\,, 0 < q, s \leq  \infty\,.
\end{equation}
\noindent
The space \(L^{p,q}(m)\) is separable if and only if \(q \not= \infty\). The 
class of functions \(f\in L^0(m)\) satisfying \(\lim_{\\ s\rightarrow 0}\
 [f^*(s)s^{\frac{1}{p}}] = 0\) is a separable closed subspace of
 \(L^{p,\infty }(m)\) which is denoted by \(L_{0}^{p,q}(m)\). We observe
here that for \(q \not=  \infty\,,\) the space \(L^{\infty , q}(m) = 0\). In the sequel 
we will follow the convention that in all discussions concerning the space \(L^{p,q}(m)\) it is 
assumed that \(q = \infty \) whenever \(p = \infty \). 

For \(w \in \mathbb {D}\,,\) and \(f \in H(\mathbb {D})\,,\) the function \(f_w\) is defined on 
\(|z| < 1/|w|\) by\(f_w(z) = f(wz)\). The space of continuous complex-valued functions on 
\(\mathbb {T}\) will be denoted by \(C(\mathbb{T})\). The function \(f_w\) is considered as 
both an analytic function on the disk \(|z| < 1/|w|\,,\) and as a function in \(C(\mathbb{T})\).
For \(0 < r < 1\,,\) the functions \(f_r\) are called the dilations of \(f\). Recall that the means
\(M_p(r, f)\) are defined in the usual way by \(M_p(r,f) = (\int_\mathbb{T} |f(r(z)|^p \, dm(z))^{1/p}\),
\(0 < p < \infty\) and \(M_\infty (r,f) = \sup_{z\in \mathbb{T}}|f(rz)| \). 
The Hardy space \(H^p\) is defined as \(H^p = \{f \in H(\mathbb {D}) : ||f||_{H^p} < \infty \}\), 
where \(||f||_{H^p} = \sup_{0 < r < 1} M_p(r, f).\) The Nevanlinna class \(N\) is the subclass
of functions \(f \in H(\mathbb {D})\) for which \(\sup_{0 < r < 1} \int _{\mathbb {T}} \log^+|f_r(z)| \,dm(z) < \infty \,.\)
Functions in \(N\) are known to have non-tangential limits \(m\mbox{-a.e. on} \mathbb{T}\).
Consequently every \(f \in H(\mathbb {D})\) determines a boundary value function which
we also denote by \(f\). Thus \(f(z) = \lim_{r\rightarrow1^-}f_r(z)\), \(m\)-a.a. \(z \in \mathbb {T}\,.\)
The Smirnov class \(N^+\) is the subclass of \(N\) consisting of those functions \(f\) for which
\(\lim_{r\rightarrow 1^-} \int _{\mathbb {T}} \log^+|f_r(z)| \, dm(z) = \int _{\mathbb {T}} \log^+ |f(z)| \, dm(z)\,.\)

It follows from standard \(H^p\)-theory that a function \(f \in H(\mathbb {D)}\) belongs to \(H^p\) 
if and only if \(f\in N^+\) with boundary value function in \(L^p(m)\,,\) in which case \(||f||_{H^p} = ||f||_p\,,\) [10].
Motivated by this characterization of \(H^p\) we define the Hardy-Lorentz space
\(H^{p,q}\), \(0 < p, q \leq  \infty \) to be the space of functions \(f \in N^+\) with boundary value
function in \(L^{p,q}(m)\) and we put \(||f||_{H^{p,q}} = ||f||_{p,q}\). Then \(\{(H^{p,q}, ||\cdot ||_{H^{p,q}})\}_{0 < p, q\leq \infty }\)
is an array of quasi-Banach spaces of analytic functions on \(\mathbb{D}\) with the standard Hardy space scale as the main 
diagonal. As with \(L^{p,q}(m)\,,\) \(H^{p,q}\) is separable if and only if \(q \not= \infty \). The functions in
\(H^{p,\infty }\) with boundary value function in \(L_0^{p,\infty }(m)\) form a closed separable 
subspace of \(H^{p,\infty }\) which is denoted by \(H_0^{p,\infty }\). Analogs of the inclusion 
relations \eqref{2.1} and \eqref{2.2} hold for the Hardy-Lorentz spaces and \(H^{p,q}\hookrightarrow H_0^{p,\infty }\) for all
\(q \not=  \infty \). The polynomials are dense in \(H_0^{p,\infty }\) and in \(H^{p,q}\,,\) \(q \not= \infty \,.\)
For \(f \in H^{p,q}\,,\) \(q \not= \infty \,,\) the dilations \(f_r\rightarrow f\) in \(H^{p,q}\) as 
\(r\rightarrow 1^- \). If \(f \in H^{p,\infty }\) then the dilations \(f_r\rightarrow f\) in \(H^{p,\infty }\) as \(r\rightarrow 1^- \)
if and only if \(f \in H_0^{p, \infty }\).  We note that a similar statement can be made for the disk algebra \(A(\mathbb{D})\) 
which is defined as the subspace of \(H^\infty \) with boundary function in \(C(\mathbb{T})\). That is, the polynomials are dense in \(A(\mathbb{D})\) 
and the dilations \(f_r\rightarrow f\) in \(H^\infty\) as \(r\rightarrow 1^- \)
if and only if \(f \in A(\mathbb{D})\), [32]. 

An important result in the theory of Lorentz spaces is the identification of these spaces
with the intermediate spaces arising in the real interpolation theory of the Lebesgue
spaces. An analytic analog of this result is given in Theorem 2.1 below. Theorem 2.1 is 
one of two interpolation theorems needed for the sequel. It was proved in [15] but omitted
the endpoint case corresponding to \(H^\infty \,.\) The complete version was proved in [29], see also [44]. 


\begin {Thm}Let \(0 < \theta <1\) and for \(j = 0, 1\) let \(0 < p_j, q_j \leq \infty\) with \(p_0 \not= p_1\).\\
\\
(i)
Set \(\frac{1}{p} = \frac{1 - \theta }{p_0} + \frac{\theta }{p_1}\). Then for every \(0 < q \leq \infty \) we have,
with equivalent quasinorms,

\[(H^{p_0,q_0},H^{p_1,q_1})_{\theta ,q} = H^{p,q}\,.\] 
\\

\noindent (ii)
Set \(\frac{1}{q} = \frac{1 - \theta }{q_0} + \frac{\theta }{q_1}\). 
Then for every \(0 < p < \infty \) we have,
with equivalent quasinorms,

\[(H^{p,q_0},H^{p,q_1})_{\theta ,q} = H^{p,q}\,.\]  

\end{Thm}

The second collection of domain spaces \(E\hookrightarrow H(\mathbb{D})\) that we define are spaces 
of mixed
`` norm " type. Before introducing these spaces, we describe the fractional calculus that we 
will be using. Let \(0 < \beta < \infty \) and suppose \(f \in H(\mathbb{D})\) with Taylor series representation

\begin{equation}\label{2.3}
f(z) = \sum_{n = 0}^{\infty } a_nz^n\,,   z \in \mathbb{D} 
\end{equation}

\noindent The fractional derivative and fractional integral of \(f\) of order \(\beta \) are the functions respectively
defined at \(z \in \mathbb{D}\) by \(f^{[\beta ]}(z) = \sum_{n=0}^{\infty} \frac{\Gamma (n + \beta + 1)}{n!}\ z^n\) and 
\(f_{[\beta ]}(z) = \sum_{n=0}^{\infty }\frac{n!}{\Gamma (n+\beta +1)}\ z^n\), where \(\Gamma \) is the gamma function. The symbols \(D^\beta \) and \(D_\beta \) stand for the associated operators
defined on \(H(\mathbb{D})\) by \(D^\beta (f) = f^{[\beta ]}\) and \(D_{\beta }(f) = f_{[\beta ]}\), \(f \in H(\mathbb{D})\). We adopt the convention that for
\(-\infty < \beta < 0\), \(f^{[\beta ]} = f_{[-\beta ]}\), \(f_{[\beta ]} = f^{[-\beta ]}\) and similarly for \(D^\beta \) and \(D_\beta \). The operators
\(D^0\) and \(D_0\) are understood to be the identity on \(H(\mathbb{D})\) and \(f^{[0]} = f_{[0]} = f\).

Suppose then that \( 0 < p, q \leq \infty , 0 < \alpha < \infty \) and let \(f \in H(\mathbb{D})\). We set 
\( ||f||_ {H(p,q,\alpha )} = (\int _0^1M_p(r,f)^q(1-r)^{q\alpha -1} \,dr)^{1/q}\,, q \not= \infty \)
and \( ||f||_{H(p,\infty ,\alpha )} = \sup_{0 < r < 1} [M_p(r,f)(1-r)^\alpha\). 
Then the weighted mixed Bergman space\(H(p,q,\alpha )\) is defined as 
\(H(p,q,\alpha ) = \{f \in H(\mathbb{D}) : ||f||_{H(p,q,\alpha )} < \infty \}\).
We also define \(H_0(p,\infty ,\alpha )\) to be the subspace of functions \(f \in H(p,\infty ,\alpha )\) satisfying \(M_p(r,f)(1-r)^\alpha\rightarrow 0\)
as \(r\rightarrow{1^-}\). In this notation, \(H(p,p,1/p)\) is the standard Bergman space
\(A^p = \{f \in H(\mathbb{D}) : \int _{\mathbb{D}} |f(z)|^p \,d\nu (z) < \infty \}\), 
where \(\nu \) is Lebesgue measure on \(\mathbb{D}\). For \( - \infty < \beta < \infty \), we set 
\(||f||_{H(p,q,\alpha,\beta)} = ||f^{[\beta ]}||_{H(p,q,\alpha )}\).
The weighted mixed Bergman-Sobolev space \(H(p,q,\alpha ,\beta )\) is then defined as
\(H(p,q,\alpha ,\beta ) = \{f \in H(\mathbb{D}): ||f||_{H(p,q,\alpha ,\beta )} < \infty \}\).
Similarly \(H_0(p,\infty ,\alpha ,\beta )\) is the subspace of \(H(p,\infty ,\alpha ,\beta )\) consisting of those functions \(f\) satisfying
\(f^{[\beta ]} \in H_0(p,\infty ,\alpha )\). The spaces
\((H(p,q,\alpha ), ||\cdot ||_{H(p,q,\alpha )})\) and 
\((H(p,q,\alpha,0), ||\cdot ||_{H(p,q,\alpha,0)})\) are of course identical and we will continue to use the former notation when \( \beta = 0\).
As with \(H^{p,q}\), \(H(p,q,\alpha ,\beta )\) is separable if and only if \(q \not= \infty \). The space \(H_0(p,\infty ,\alpha ,\beta )\) is a 
closed separable subspace of \(H(p,\infty ,\alpha ,\beta )\). The polynomials are dense in \(H(p,q,\alpha ,\beta ), q \not= \infty \) and in \(H_0(p,\infty ,\alpha ,\beta )\).
If \(f \in H(p,q,\alpha ,\beta ), q \not= \infty \), then \(f_r\rightarrow f\) in \(H(p,q,\alpha ,\beta )\) as \(r\rightarrow 1^-\).
On the other hand, for \(f \in H(p,\infty ,\alpha ,\beta )\), \(f_r\rightarrow f\) in \(H(p,\infty ,\alpha ,\beta )\) as \(r\rightarrow 1^-\) if and only if
\(f \in H_0(p,\infty ,\alpha ,\beta )\), see [28], [50]. The spaces \(H(p,q,\alpha ,\beta )\) are often called mixed norm spaces. In the sequel we will simply say
\(H(p,q,\alpha ,\beta )\) is a Bergman-Sobolev spaces and \(H(p,q,\alpha )\) is a Bergman space. Many authors use an equivalent 
definition of these spaces obtained by replacing the fractional calculus operators \(D^\beta \) and \(D_\beta \,, 0 \leq \beta < \infty \), in the 
definition of \(H(p,q,\alpha ,\beta )\) with the multiplier operators \(J^\beta \) and \(J_\beta \) defined at a function \(f \in  H(\mathbb{D})\) with Taylor 
series representation \eqref{2.3} by
\(J^\beta (f)(z) = \sum_{n = 0}^{\infty } (n+1)^\beta a_n\,z^n \mbox{ and } J_\beta (f)(z) = \sum_{n = 0}^{\infty } (n+1)^{-\beta } a_n\,z^n\,,z \in \mathbb{D}\).  
Then the mixed norm space obtained using \(J^\beta \) or \(J_\beta \) is identical to the space \(H(p,q,\alpha ,\beta )\) as previously defined.
For the equivalence of the fractional calculus operators with the multiplier operators in defining the spaces \(H(p,q,\alpha ,\beta )\) as well
as proofs of the following results, the reader is referred to [16], [48], and [49].
\begin{Lem} Let \(0 < p\,, q \leq \infty \), \(0 < \alpha , \beta < \infty \). Then the following mappings are continuous surjective isomorphisms.\\

\noindent
(i) \(D^\beta : H(p,q,\alpha )\rightarrow H(p,q,\alpha + \beta ),\)\\
(ii) \(D^\beta : H_0(p,\infty ,\alpha )\rightarrow H_0(p,\infty ,\alpha + \beta ),\)\\
(iii) \(D_\beta : H(p,q,\alpha )\rightarrow H(p,q,\alpha - \beta ) \mbox { for } \beta < \alpha\),\\
(iv) \(D_\beta : H_0(p,\infty ,\alpha )\rightarrow H_0(p,\infty ,\alpha - \beta ) \mbox { for } \beta < \alpha\).
\end{Lem}

Lemma 2.1 and the definition of \(H(p,q,\alpha ,\beta )\) imply the following.
\begin{Lem} Let \(0 < p, q \leq \infty \,,\) \(0 < \alpha < \infty, - \infty < \beta < \infty  \). Then for \( -\infty < \gamma < \alpha \,,\) the 
following identifications hold with equivalent quasinorms.\\

\noindent
(i) \(H(p,q,\alpha ,\beta ) = H(p,q, \alpha - \gamma , \beta - \gamma )\,,\)\\
(ii) \(H_0(p,\infty ,\alpha ,\beta ) = H_0(p,\infty ,\alpha  - \gamma ,\beta -\gamma )\,.\) 
\end{Lem}

In particular, Lemma 2.2 implies \(H(p,q,\alpha ,\beta) = H(p,q,\alpha -\beta )\) and \(H_0(p,\infty ,\alpha ,\beta ) = H_0(p,\infty ,\alpha -\beta )\)
if \(-\infty < \beta < \alpha\). Lemma 2.3 which follows represents an extension of Lemma 2.1 to the spaces \(H(p,q,\alpha ,\beta)\).  
\begin{Lem} Let \(0 < p\,, q \leq \infty \), \(0 < \alpha < \infty, - \infty < \beta, \gamma  < \infty  \). Then the following mappings are
continuous surjective isomorphisms.\\

\noindent
(i) \(D^\gamma : H(p,q,\alpha ,\beta )\rightarrow H(p,q,\alpha + \gamma,\beta )\) for \(\gamma > - \alpha, \)\\
(ii) \(D^\gamma : H_0(p,\infty ,\alpha ,\beta )\rightarrow H_0(p,\infty ,\alpha + \gamma,\beta )\) for \(\gamma > - \alpha\),\\
(iii) \(D^\gamma : H(p,q,\alpha ,\beta )\rightarrow H(p,q,\alpha, \beta  - \gamma),\)\\
(iv) \(D^\gamma : H_0(p,\infty ,\alpha ,\beta )\rightarrow H_0(p,\infty ,\alpha, \beta  - \gamma).\)
    
\end{Lem}  

The second interpolation result we need is for the Bergman-Sobolev spaces and is due to Fabrega and Ortega, see [14].
\begin{Thm}Let \(0 < p < \infty\,,0< q_j \leq \infty\,, 0 < \alpha _j, \beta < \infty \,,
j=0,1 \mbox{ and suppose } \alpha _0 \not= \alpha _1\). 
Let \(0 < \theta < 1 \,, 0 < q \leq \infty\) and set \(\alpha = (1 - \theta )\alpha _0 + \theta \alpha _1\). Then we have, with equivalent quasinorms,

\[(H(p,q_0,\alpha _0,\beta ), H(p,q_1,\alpha _1,\beta ))_{\theta, q}= H(p,q,\alpha ,\beta )\,.\] 
\end{Thm}
We also need the following embedding theorems. The first of these is due to Flett [16] and indicates how Bergman-Sobolev spaces embed 
in the standard Hardy spaces. The second result is a well-known embedding theorem of Hardy and Littlewood, see [10].
\begin{Thm} Let \(0 < p < s < \infty\), \(0 < q \leq s\), \(\beta > 1/p - 1/s\). Then
\[H(p,q,\beta + 1/s -1/p, \beta )\hookrightarrow\ H^s\,.\]
\end{Thm}
\begin{Thm} 
Let \(0 < p < s \leq \infty \), \(p \leq t \leq \infty \). Then
\[H^p\hookrightarrow\ H(s,t,1/p- 1/s)\,.\]  
\end{Thm}

\setcounter{equation}{0}
\section{Multipliers}
Let \(\mathbb{N}_0\) denote the set of nonnegative integers and let \(W\) denote the space of complex sequences indexed by  \(\mathbb{N}_0\).
We always consider \(W\) as being equipped with the topology of pointwise convergence. With this topology, \(W\) is a Frechet space. A topological 
vector space \(X\) satisfying \(X\hookrightarrow W\) is called a K-space. An FK-space is a K-space which is also an F-space. In particular, spaces which
are both K-spaces and Frechet spaces are FK-spaces. \(W\) is also a topological algebra under the natural product of coordinate-wise
multiplication. Thus, for \(w = \{w_n\}, \lambda =\{\lambda _n\}\), the product \(\lambda w\) is defined by \(\lambda w = \{\lambda _n w_n\}\).
It will sometimes be convenient to use the symbol \(B\) for the product map so that \(B(\lambda ,w) = \lambda w\). Then 
\(B: W \times W\rightarrow W\) is a continuous bilinear operator. For fixed \(\lambda \in W\), we will write \(B_\lambda \) for the continuous 
linear operator \(B_\lambda : W\rightarrow W\) defined by \(B_\lambda (w) = \lambda w\), \(w \in W\). Suppose now that \(E\) and \(X\)
are a pair of vector subspaces of \(W\). An element \(\lambda \in W\) is said to be a multiplier of \(E\) into \(X\) if \(\lambda w \in X\) for
every \(w \in E\). The set of multipliers from \(E\) into \(X\) is denoted by either of the symbols \((E,X)\) or \(E^X\). Thus \(\lambda \in (E,X)\)
if and only if the linear operator \(B_\lambda \) maps \(E\) into \(X\). (Consequently, the bilinearity of \(B\) gives \((E,X)=\bigcap_{w\in E}(B_w^{-1}(X)\)).
If \(E\) and \(X\) are FK-spaces, an argument based on the Closed Graph Theorem shows that \((E,X)\), or more precisely \(\{B_\lambda : \lambda \in (E,X)\}\), 
is a subspace of \( \cal L\)\((E,X)\), the space of continuous \(X\)-valued linear operators on \(E\). The space \((E,X)\) is sometimes called the \(X\)-dual 
of E. The second \(X\)-dual of \(E\) is the space \(E^{XX} = (E^X)^X\). If   \(E^{XX} = E\) then E is said to be \(X\)-reflexive or \(X\)-perfect. We record 
some of the basic properties of multiplier spaces in the form of a lemma. We omit the obvious proof.  
\begin{Lem} Let \(A\), \(B\), \(C\), \(E\) be vector subspaces of \(W\) with \(A\subset B\). Then\\

\noindent
(i) \(B^C\subset A^C\,,\)\\
(ii) \(C^A\subset C^B\,,\)\\
(iii) \((A,C)\subset (C^E, A^E)\,.\)
\end{Lem}

For quasi-Banach spaces \((E,||\cdot ||_E)\), \((X,||\cdot ||_X)\) \(\hookrightarrow W\), the operator quasinorm is defined at an operator 
\(L \in \cal L\)\((E,X)\) in the standard way by \(||L||_{{\cal L }(E,X)}\)\(=\sup\{||L(w)||_X : w \in E\,,||w||_E \leq 1\}\). Then
\((\cal L\)\((E,X)\,,\) \(||\cdot ||_{{\cal L }(E,X)})\) is a quasi-Banach space containing \((E,X)\) as a closed subspace. In particular, \((E,X)\) is a quasi-Banach space under the quasinorm \(||\cdot ||_{(E,X)}\)
defined at \(\lambda \in (E,X)\) by \(||\lambda ||_{(E,X)} = ||B_\lambda ||_{{\cal L} (E,X)}\).

We regard \(H(\mathbb{D})\) as a subspace of \(W\) by identifying functions in \(H(\mathbb{D})\) with their Taylor coefficient sequences. 
Thus, a function \(f \in H(\mathbb{D})\) with Taylor series representation \eqref{2.3} is identified with the sequence \(a=\{a_n\}\).
\(H(\mathbb{D})\) is a Frechet space when equipped with the topology of uniform convergence on compact subsets of \(\mathbb{D}\).
In addition, we note that \(H(\mathbb{D})\) is a K-space and hence any F-space \(E\) satisfying \(E\hookrightarrow  H(\mathbb{D})\) is a
FK-space. In [55] it is shown that the product map \(B\) on \(W\times W\) restricts to a bilinear operator \(H(\mathbb{D})\times H(\mathbb{D})\rightarrow H(\mathbb{D})\).
The symbol \(c\) will denote the Cauchy function defined by \(c(z) = (1-z)^{-1}\,,z \in \mathbb{D}\). Since \(c \in H(\mathbb{D})\,,\)
it follows that \((H(\mathbb{D}),H(\mathbb{D})) = H(\mathbb{D})\). If \(f\), \(g\) \(\in H(\mathbb{D})\), then \(B(f,g)\) is commonly denoted by \(f\ast g\) and is 
called the Hadamard product of \(f\) and \(g\). Thus, if \(f\), \(g\) \(\in H(\mathbb{D})\), with Taylor series representations \(f(z) = \sum_{n = 0}^{\infty} a_n z^n\)
and \(g(z) = \sum_{n = 0}^{\infty} b_n z^n\), \(z \in \mathbb{D}\), then \(B(f,g) = f\ast g \in H(\mathbb{D})\) has 
Taylor series representation \(B(f,g)(z) = (f\ast g)(z) = \sum_{n = 0}^{\infty } a_n b_n z^n\,, z \in \mathbb{D} \) and, as a sequence, \(B(f,g) = f\ast g = \{a_n b_n\}\).\\
We introduce some notation. For \(n \in \mathbb{N}_0\) and \(z \in \mathbb{D}\) we set \(u_n(z) = z^n\). Let \(f \in H(\mathbb{D})\)
with Taylor series representation given by \eqref{2.3}. For \(N \in \mathbb{N}_0\) we write \(S_N(f)\) for the partial sum function
\(S_N(f)(z) = \sum_{n = 0}^{N} a_n z^n\,, z \in \mathbb{D}\,.\)
\noindent
In the sequel we will be interested in the multiplier spaces \((E,X)\) where E is a Hardy-Lorentz
space and \(X\) is a FK-space, \(X\hookrightarrow H(\mathbb{D})\). At this point we would like to consider \((E,X)\) for
some specific target spaces \(X\) and for a fairly general class of domain spaces \(E\). One of our choices for
\(X\) is the space \(AS(\mathbb{D})\) of Abel summable sequences. Recall that the element \(w = \{w_n\} \in W\)
is said to be Abel summable if \(\lim_{r\rightarrow 1^-} \sum_{n = 0}^{\infty } w_n r^n \) exists. The space  \(AS(\mathbb{D})\) 
is a Frechet space [45], with respect to the topology induced by the family \(\{\rho _n : n \in \mathbb{N}_0 \mbox{ or } n = \infty \}\)  
of seminorms, where \(\{r_n\}\) is a fixed sequence in \((0,1)\) strictly increasing to \(\infty\), and for \(n \in \mathbb{N}_0\) and \(w = \{w_n\} \in AS(\mathbb{D})\), 
\(\rho _\infty(w) = \sup_{0 < r < 1} | \sum_{k = 0}^{\infty } w_k r^k | \mbox{ and } \rho _n(w) = \sum_{k = 0}^{\infty } |w_k| r_n^k\,.\)    
Furthermore, we have  \(AS(\mathbb{D})\hookrightarrow H(\mathbb{D})\) and hence \(AS(\mathbb{D})\) is a FK-space. The \(AS(\mathbb{D})\)-dual
of a vector subspace \(E\) of \(W\) is known as the Abel dual of \(E\) and will be denoted by \(E^a\). The second
Abel dual of \(E\) is the space \((E^a)^a\) and will be denoted by \(E^{aa}\). If \(E^{aa} = E\), then \(E\) is said to be Abel reflexive.
Note that the functional \(\Psi \) on  \(AS(\mathbb{D})\) defined by 

\begin{equation}\label{3.1}
\Psi (a) = \lim_{r\rightarrow 1^-} \sum_{n = 0}^{\infty } a_n r^n\,, a = \{a_n\} \in AS(\mathbb{D})
\end{equation}

\noindent belongs to \(AS(\mathbb{D})^*\,.\)\\
\noindent 
Propositions 3.2 through 3.4 below describe the relationship between the spaces \(E^*\), \((E,H^\infty )\), and \(E^a\) for a certain general type of
FK-space \(E\). First we need the following.
\begin{Prop} Suppose that \(E\) is a FK-space satisfying\\

\noindent
(i) \(E\hookrightarrow H(\mathbb{D})\,,\)\\
(ii) \(u_n \in E\mbox{ for all } n \in \mathbb{N}_0\,,\)\\
(iii) \(c_w \in E\mbox{ for all } w \in \mathbb{D}\,,\)\\
(iv) For all \(w \in \mathbb{D}\,,\) \(\{S_N(c_w)\}\mbox{ converges to } c_w\mbox { in } E\mbox { as } N\rightarrow \infty\,.\)\\

\noindent
Then for any topological vector space \(X\), operator \(T \in \cal L\)\((E,X)\), and \(w \in \mathbb{D}\), the \(X\)-valued series

\begin{equation}\label{3.2}\sum_{n = 0}^{\infty } x_n w^n\,, x_n = T(u_n)\,, n \in \mathbb{N}_0\,,
\end{equation}  

\noindent
converges in \(X\) to \(T(c_w)\).
\end{Prop}

\noindent
\begin{proof}
The partial sums of the series \eqref{3.2} satisfy
\begin{equation}\label{3.3}
\sum_{n = 0}^{N} x_n w^n = T(S_N(c_w))\,.
\end{equation} 
Since T is continuous the lemma follows from (iv) and \eqref{3.3}.\\
\end{proof}

\indent
Suppose then that \(E\) is a FK-space satisfying (i) through (iv) of Proposition 3.1, that \(X\) is a topological vector space and 
\(T \in \cal L\)\((E,X)\). Define a function \(g_T\) on \(\mathbb{D}\) by \(g_T(w) = T(c_w)\,,w\in \mathbb{D}\). 
Then \(g_T\) is a well-defined \(X\)-valued function on \(\mathbb{D}\) with power series representation given by \eqref{3.2}.
The function \(g_T\) is called the analytic or Cauchy transform of \(T\). The linear operator \(T\rightarrow g_T\,,\) taking \(\cal L\)(\(E,X)\) 
into the space of \(X\)-valued power series on \(\mathbb{D}\,,\) will also be referred to as the analytic or Cauchy transform on  \(\cal L\)(\(E,X)\).
Consider now the case when \(X\) is a FK-space satisfying \(X\hookrightarrow H(\mathbb{D})\) and \(T = B_\lambda \) for some multiplier \(\lambda 
\in (E,X)\). For this case, if \(\lambda = \{\lambda _n \}\) as a sequence of complex numbers then \(T(u_n) = B_\lambda (u_n) = \lambda u_n\).
Since \(\lambda u_n\) is identified with the sequence having \(\lambda _n\) in the \(n\)-th entry and 0 elsewhere, we will write \(\lambda _n\) in
place of \(T(u_n)\,,\) \(n \in \mathbb{N}_0\,,\) and \(g_\lambda\) in place of \(g_T\). The other situation we will be considering is when
\(X\) is the complex field and \(T = \varphi \in E^*\). For this case we may form the sequence \(\lambda =\{\lambda _n\}\,,\) where
\(\lambda _n = \varphi(u_n)\,, n \in \mathbb{N}_0\) and we again write \(g_\lambda \) in place of \(g_T\). Proposition 3.1 ensures that \(g_\lambda \in
H(\mathbb{D})\). But, as was previously noted, \(H(\mathbb{D}) = (H(\mathbb{D}),H(\mathbb{D}))\). Since
\((H(\mathbb{D}),H(\mathbb{D}))\subset (E,H(\mathbb{D}))\) by Lemma 3.1(i), it follows that \(\varphi \in E^*\)
induces the multiplier \(\lambda \in (E,H(\mathbb{D}))\) with analytic transform \(g_\lambda \). In fact, by convention, we
have the identification \(\lambda \leftrightarrow g_\lambda \). Let us note here that it is possible to have \(\lambda =0\) for \(\varphi \not= 0\,,\)
so that in general the analytic transform on \(\cal L\)\((E,X)\) is not one-to-one. However we do have the following. 
\begin{Prop} Suppose that in addition to satisfying conditions (i) through (iv) of Proposition 3.1, the FK-space \(E\) satisfies \\

\noindent
(i) \(f_w \in E\mbox{ for each }w \in \mathbb{D}\mbox{ and } f \in E\,,\)\\
(ii) for each \(f \in E\,,\mbox{ the set } \{f_w : w \in \mathbb{D}\}\mbox{ is bounded in } E\,,\)\\
(iii) \(\{S_N(f_w)\}\mbox{ converges to } f_w \mbox{ in } E\mbox{ for every } w \in \mathbb{D}\mbox{ and } f \in E\,.\)\\

\noindent
Let \(\varphi \in E^*\) with induced multiplier \(\lambda = \{\lambda _n\}\,,\lambda _n = \varphi (u_n)\) and analytic transform \(g_\lambda \).
Then for each \(w \in \mathbb{D}\) and \(f \in E\),

\begin{equation}\label{3.4}
\varphi(f_w) = (f\ast g_\lambda )(w)\,.
\end{equation}

\noindent
Consequently, \(\lambda \in (E,H^\infty )\). Furthermore, the analytic transform on \(E^*\) is one-to-one whenever
\(E\) satisfies the additional property \\

\noindent
(iv) for every \(f \in E\), the dilations \(f_r\), \(0 < r < 1\), converge to \(f\) in \(E\) as \( r\rightarrow 1^-\,.\)

\end{Prop}
\noindent
\begin{proof}
Let \(\varphi \in E^*\) and put \(\lambda = \{\lambda _n\}\,, n \in \mathbb{N}_0\). Then , by the comments made in the paragraph
following the proof of Proposition 3.1, we have \(g_\varphi = g_\lambda \in (E, H(\mathbb{D}))\). If \(w \in \mathbb{D}\) and \(f \in E\) has
Taylor series representation \eqref{2.3} then the continuity of \(\varphi\) and conditions
(i) and (iii) of Proposition 3.2 yield
\(\varphi(f_w) = \varphi(\lim_{n\rightarrow \infty }S_N(f_w)) = \lim_{N\rightarrow\infty }\varphi(S_N(f_w) = 
\lim_{N\rightarrow\infty }\varphi(\sum_{n=0}^{N} a_n w^n u_n) = \lim_{N\rightarrow \infty }\sum_{n=0}^{N} a_n\lambda _n w^n
= \lim_{N\rightarrow \infty }S_N(f\ast g_\lambda )(w) = f\ast g_\lambda (w)\) which is \eqref{3.4}. Then \eqref{3.4}, the continuity of \(\varphi\), and condition (ii) imply \(\lambda \in (E,H^\infty )\).
Finally, suppose \(E\) satisfies (iv). It then follows from this property and \eqref{3.4}, that if \(\lambda \) is the
induced multiplier for \(\varphi_j \in E^*\,, j = 1,2\) we have \(\varphi_j(f) = \varphi_j(\lim_{r \rightarrow 1^-}f_r) 
= \lim_{r\rightarrow 1^-}\varphi_j(f_r) = \lim_{r\rightarrow 1^-}(f\ast g_\lambda )(r)\), so that \(\varphi_1 = \varphi_2\).
\end{proof}

In view of the last two propositions, we may regard \(E^*\subset(E,H^\infty )\) for any FK-space
satisfying (i) through (iv) of Propositions 3.1 and 3.2. Note that if \(E\) is a FK-space with  \(A(\mathbb{D})\hookrightarrow E\) then \(E\) satisfies conditions (ii) 
through (iv) of Proposition 3.1 and conditions (i) through (iii) of Proposition 3.2. Also in [42] it is shown
that for the A-spaces, which form a large class of quasi-Banach spaces \(E\hookrightarrow H(\mathbb{D})\),
condition (iv) of Proposition 3.2 is equivalent to the density of the polynomials in E.
\begin{Prop} Suppose \(E\) is an FK-space. Then \(E^a\subset E^*\). 
\end{Prop}
\begin{proof}
Let \(\lambda \in E^a\) and put \(\varphi_\lambda  = \Psi \circ B_\lambda \), where \(\Psi \) is
the functional \eqref{3.1}. The mapping \(\lambda \rightarrow\varphi_\lambda \) is a one-to-one linear operator 
taking \(E^a\) into \(E^*\) and we may identify \(E^a\) with its image in \(E^*\).
\end{proof}
\begin{Prop}Suppose \(E\) is an FK-space satisfying conditions (i) through (iv) of Propositions 3.1 and 3.2.
Then \((E,H^\infty ) = (E,A(\mathbb{D})) \subset E^a \).
\end{Prop}
\begin{proof}

Let \(\lambda = \{\lambda _n\}\in (E,H^\infty )\). Let \(f\in E\) have Taylor series representation
\eqref {2.3}. By Proposition 3.2(iv) and the continuity of \(B_\lambda\) we have \(B_\lambda(f_r) \rightarrow B_\lambda(f)\) in 
\(H^\infty \). Hence \(B_\lambda (f) \in A(\mathbb{D})\). Therefore, \((E,H^\infty ) \subset (E,A(\mathbb{D}))\). By Lemma 3.1(ii) the reverse inclusion holds.
Thus  \((E,H^\infty ) = (E,A(\mathbb{D}))\).  Finally, the continuity of \(B_\lambda (f)\) at 1
 gives \(\lim_{r\rightarrow 1^-} \sum_{n=0}^{\infty }a_n\lambda _n r^n = \lim_{r\rightarrow 1^-} B_\lambda (f)(r) = B_\lambda (f)(1)\). Hence \(\lambda \in E^a\).
\end{proof} 
\begin{Cor} Let \(E\) be a FK-space satisfying conditions (i) through (iv) of Propositions 3.1 and 3.2.
Then \(E^* = E^a = (E,H^\infty )\).

\end{Cor}
\begin{Cor} Let \(0 < p, q < \infty \). Then 

\noindent
(i) \((H^{p,q})^* = (H^{p,q})^a = (H^{p,q},H^\infty )\,,\)\\
(ii) \((H_0^{p,\infty })^* = (H_0^{p,\infty })^a = (H_0^{p,\infty },H^\infty )\,.\)
\end{Cor}

\setcounter{equation}{0}
\section{A Nested Embedding Theorem For Hardy-Lorentz Spaces}
Our main tool for identifying certain multiplier spaces \((H^{p,q},X)\) is the following nested embedding theorem for \(H^{p,q}\).
Its proof consists of using Theorems 2.1 and 2.2 to interpolate Theorems 2.3 and 2.4.
\begin{Thm}
Let \(0 < p_0 < p < s \leq \infty\,,\, 0 < q \leq t \leq \infty\) and \(\beta > 1/p_0 - 1/p\). Then\\

\noindent
(i) \(H(p_0,q,\beta + 1/p-1/p_0,\beta )\hookrightarrow H^{p,q}\hookrightarrow H(s,t,1/p-1/s)\,.\)\\
(ii) \(H_0(p_0,\infty ,\beta + 1/p - 1/p_0,\beta )\hookrightarrow H_0^{p,\infty }\hookrightarrow H_0(s,\infty ,1/p-1/s)\,.\)  
\end{Thm}
\begin{proof}
(i) Choose \(s_j\), \(j = 0, 1\), such that \(0 < p_0 < s_0 < p < s_1 < \infty\). Let us further stipulate that for the case \(\beta 
\leq 1/p_0\) we require that \(0 < s_1 < \frac{p_0}{1-\beta p_0}\). This ensures that \(\beta > 1/p_0 - 1/s_1\). We can then apply Theorem 2.4 
to obtain embeddings
\begin{equation}\label{4.1}
H(p_0,s_j,\beta + 1/s_j - 1/p_0,\beta )\hookrightarrow H^{s_j}\,,j = 0, 1
\end{equation} 
\noindent
Interpolation of \eqref{4.1} results in the embeddings
\begin{equation}\label{4.2}
\Big(H(p_0,s_0,\beta + 1/s_0 - 1/p_0, \beta ),H(p_0,s_1,\beta + 1/s_1 - 1/p_0)\Big)_{\theta,q}\hookrightarrow (H^{s_0},H^{s_1})_{\theta ,q}\,,
\end{equation}
\\
\noindent 
for all \(0 < \theta < 1, 0 < q \leq \infty\). Then the first embedding in (i) follows by choosing \(\theta \) to satisfy
\(\frac{1}{p} = \frac{1 - \theta}{s_0} + \frac{ \theta}{s_1}\)
and applying Theorems 2.1 and 2.2 to \eqref{4.2}. The proof of the second embedding in (i) is similar and is in [30]. We omit the details.\\
\noindent(ii) Again, we prove only the first embedding in (ii) since the proof of the second embedding in (ii) is similar and is also in [30].
Thus, let\\
\(f\in H_0(p_0,\infty ,\beta + 1/p - 1/p_0, \beta )\). Then \(f \in H^{p,\infty }\) by Theorem 4.1(i). In order to show \(f \in H_0^{p,\infty }\), it is enough to show
the dilations \(f_r\) converge to \(f\) in \(H^{p,\infty }\) as \(r\rightarrow1^-\). But the functions \(f_r\) converge to \(f\) in \(H(p_0,\infty ,\beta + 1/p - 1/p_0,\beta )\)
and this fact combined with (i) implies \(f_r\rightarrow f\) in \(H^{p,\infty }\). Hence \(f \in H_0^{p,\infty }\).
\end{proof}

Recall that if \(E\) is a quasi-Banach space with separating dual \(E^*\), then there exists a unique Banach space \(Y\)
in which \(E\) embeds as a dense subspace and for which \(Y^* = E^*\). The space \(Y\) is called the Banach envelope
of \(E\) and is denoted by \([E]_1\). In [33] we identified the Banach envelopes and dual spaces of the spaces \(H^{p,q}\) and \(H_0^{p,\infty }\) 
for indices in the range \(0 < p < 1\), \(0 < q < \infty \). The specific result was the following. 
\begin{Thm}Let \(0 < p < 1 \), \(0 < q < \infty\). Set \(q_{\ast }=max(1,q)\) and let \(q'\) be the H\"{o}lder conjugate of \(q_{\ast }\), \(1/q_{\ast }+1/q' = 1\).
Then\\

\noindent (i)\([H^{p,q}]_1 = H(1,q_{\ast },1/p - 1)\mbox{ and } (H^{p,q})^* = H(\infty,q',1,1/p)\,,\)\\ 
(ii) \([H_0^{p,\infty }]_1 = H_0(1,\infty ,1,1/p -1)\) and \((H_0^{p,\infty })^* = H(\infty ,1,1,1/p)\). 
\end{Thm}

\noindent Note that for \(p = q\), (i) becomes

\begin{equation}\label{4.3}
[H^p]_1 = H(1,1,1/p-1)\mbox{ and } (H^p)^* = H(\infty ,\infty ,1,1/p) 
\end{equation}

\noindent This is the well-known Duren-Romberg-Shields Theorem [11]. The proof of Theorem 4.2(i) given in [33] essentially
consisted of two steps. The first step was the establishment of the second embedding in Theorem 4.1(i). The 
second step was a constructive proof of the embedding\\
\((H^{p,q})^*\hookrightarrow H(\infty ,q',1,1/p)\). The proof
of Theorem 4.2(ii) in [33] was carried out in an analogous fashion. A short proof of Theorem 4.2 can be based on Theorem 4.1
and Theorem 4.3 below. Theorem 4.3 is a general Banach envelope-duality theorem for separable
Bergman-Sobolev spaces and is due to the efforts of several authors. For statements and proofs of Theorem 4.3
for Bergman spaces in some specific cases the reader is referred to [1], [2], [3], [6], [7], [9], [11], [17], [18], [19], [23],
[24], [25], [35], [37], [39], [46], [47], [50], [52], [54], [56]. Pavlovic's paper [42] contains a very general and complete version 
of Theorem 4.3 for the case \(\beta = 0\). To obtain the result for Bergman-Sobolev spaces one uses the validity of Theorem 4.3 for Bergman spaces with Lemmas 2.1
through 2.3.
\begin{Thm} Let \(0 < p \leq \infty , 0 < q, \alpha < \infty , - \infty < \beta < \infty \). Let \(p_0 = min(1,p), p_1 = max(1,p)\mbox{ and } q_1 = max(1,q)\). 
Let \(1/p_1 + 1/p_1' = 1/q_1 + 1/q_1' =1\). Then\\

\noindent (i) \([H(p,q,\alpha, \beta)]_1 = H(p_1,q_1,\alpha + 1/p_0 -1, \beta)\,,\) \\
(ii) \((H(p,q,\alpha, \beta))^* = H(p_1',q_1',1,\alpha - \beta + 1/p_0)\,,\)\\
(iii) \([H_0(p,\infty ,\alpha ,\beta )]_1 = H_0(p_1,\infty ,\alpha + 1/p_0 -1,\beta )\,,\)\\
(iv) \((H_0(p,\infty ,\alpha ,\beta)^* = H(p_1',1,1,\alpha - \beta + 1/p_0)\,.\)
\end{Thm}  

For \(0 < p \leq \infty , 0 < q, \alpha < \infty , - \infty < \beta < \infty \,,\) the spaces \(H(p,q,\alpha ,\beta )\)
and \(H_0(p,\infty ,\alpha ,\beta )\) satisfy conditions (i) through (iv) of Propositions 3.1 and 3.2. So by Corollary 3.1,
\(H(p,q,\alpha ,\beta )^* = H(p,q,\alpha ,\beta )^a = (H(p,q,\alpha ,\beta ), H^\infty )\) and similarly for
\(H_0(p,\infty ,\alpha ,\beta )\). Thus duality and Abel duality coincide for the separable Bergman-Sobolev spaces.
More explicitly, say in the case of Theorem 4.3(ii), if \(\Lambda \in H(p,q,\alpha ,\beta )^*\) and 
\(f \in H(p,q,\alpha ,\beta )\) has Taylor series representation \eqref{2.3}, then the proof of Theorem 4.3(ii)
shows that 

\begin{equation}\label{4.4}
\Lambda (f) = \lim_{r\rightarrow1^-}\sum_{n = 0}^{\infty } a_n \lambda _n r^n\,, \lambda _n = \Lambda (u_n)\,, n \in \mathbb{N}_0.
\end{equation}
  
\noindent Furthermore, the analytic transform \(g_\lambda \) of the sequence \(\lambda = \{\lambda _n\}\) in \eqref{4.4}
satisfies \(g_\lambda \in H(p_1,1,1,\alpha - \beta + 1/p_0)\,,||g_\lambda ||_{H(p_1,q_1',1,\alpha - \beta + 1/p_0)} \sim ||\Lambda ||_{H(p,q,\alpha ,\beta )^*}\).
Conversely, if \(g\in H(p_1,1,1,\alpha - \beta + 1/p_0)\) has Taylor coefficient sequence \(\lambda = \{\lambda _n\}\),
then we may define \(\Lambda _g\) as in \eqref{4.4}. The resulting functional \(\Lambda _g\) belongs to
\(H(p,q,\alpha ,\beta )^*\), has analytic transform \(g\), and satisfies\\
\(||\Lambda _g||_{H(p,q,\alpha ,\beta )^*}\sim ||g||_
{H(p_1,1,1,\alpha - \beta + 1/p_0)}\). Theorem 4.3 also implies that the spaces \(H(p,q,\alpha ,\beta )\,, 1 \leq p \leq \infty , 1 < q < \infty\)
are reflexive with the properties of reflexivity and Abel reflexivity being the same for these spaces.

To see how Theorem 4.2 follows from Theorems 4.1 and 4.3, let \(0 < p < 1\), choose \(0 < p_0 < p\) and take \(s = 1\,, t =\mbox {max}(1,q) = q_\ast \) 
in Theorem 4.1(i) to obtain the nested embedding

\begin{equation}\label{4.5} 
H(p_0,q,\beta + 1/p - 1/p_0, \beta )\hookrightarrow H^{p,q}\hookrightarrow H(1,q_\ast ,1/p -1)
\end{equation}

\noindent Applying the functor \([\cdot ]_1\) to \eqref {4.5} we find
\[[H(p_0,q,\beta + 1/p - 1/p_0,\beta )]_1\hookrightarrow [H^{p,q}]_1\hookrightarrow [H(1,q_\ast ,1/p -1)]_1\,.\]
Since \(H(1,q_\ast ,1/p -1)\) is a Banach space,

\begin{equation}\label{4.6}
[H(1,q_\ast ,1/p -1)]_1 = H(1,q_\ast ,1/p -1)\,.
\end{equation}

\noindent Using Theorem 4.3, we also find
\begin{equation}\label{4.7}
[H(p_,q,\beta + 1/p - 1/p_0,\beta )]_1 = H(1,q_\ast ,\beta + 1/p -1,\beta )\,.
\end{equation}
\noindent
But the spaces on the right-hand sides of \eqref{4.6} and \eqref{4.7} are identical by Lemma 2.2(i).
This establishes the first equality in Theorem 4.2(i) and hence the second inequality as well via Theorem 4.3(ii).
The proof of Theorem 4.2(ii) is similar.\\

Combining Theorems 4.2 and 4.3 one sees that the Banach envelopes of the spaces \(H^{p,q}\) are reflexive
and Abel reflexive for \(0 < p < 1 < q < \infty \). Similarly, from Corollary 3.2 and Theorem 4.2, we deduce that 
for \(0 < p < 1\,,\)

\begin{equation}\label{4.8}
(H_0^{p,\infty })^a = (H_0^{p,\infty })^* = H(\infty ,1,1,1/p)\,.
\end{equation}

\noindent It then follows from \eqref{4.8}, Theorem 4.3(ii), and Lemma 2.2(i),that

\begin{equation}\label{4.9}
(H_0^{p,\infty })^{aa} = H(1,\infty ,1/p -1)
\end{equation}

\noindent That we also have \((H^{p,\infty })^a = H(\infty ,1,1,1/p)\) is a consequence of the following result of Shi [51].
\begin{Lem} Let \(0 < \alpha < \infty \). Then\\

\noindent(i) \(H(1,\infty ,\alpha )^a = H_0(1,\infty ,\alpha )^a\,,\)\\
(ii) \(H(1,\infty ,\alpha )^{aa} = H(1,\infty ,\alpha )\,.\)\\
\end{Lem}

\begin{Cor}Let \(0 < p < 1\). Then\\

\noindent(i) \((H^{p,\infty })^a = H(\infty,1,1,1/p)\,,\)\\
(ii) \((H^{p,\infty })^{aa} = H(1,\infty ,1/p -1)\,.\)\\  
\end{Cor}
\begin{proof}
(i) Using Lemma 3.1(i) twice, \eqref{4.8}, Theorems 4.2, 4.3 and Lemma 4.1(i) we obtain
\begin{eqnarray}\label{4.10}
H(1,\infty ,1/p -1)^a  \subset  (H^{p,\infty })^a \subset (H_0^{p,\infty })^a = (H_0^{p,\infty })^* & = & H_0(1,\infty ,1/p -1)^*\notag\\
                                                                                                                                               & = & H_0(1,\infty ,1/p -1)^a\notag\\
                                                                                                                                               & = & H(1,\infty ,1/p -1)^a \notag\\ 
\end{eqnarray}
Since the endpoint spaces in \eqref{4.10} are the same, (i) follows from \eqref{4.8} and \eqref{4.10}.\\
\noindent (ii) This follows from \eqref{4.9} and \eqref{4.10}.
\end{proof}

\setcounter{equation}{0}
\section{Multipliers of \(H^{p,q}\mbox { and } H_0^{p,\infty}\mbox{ into }\ell^s\,,\\ 0 < p < 1\,,\)\(0 < q < \infty\,,0 < s \leq \infty\,.\)}

In this section we determine the multiplier spaces \((E,X)\) where \(E\) is either \(H^{p,q}\) or \(H_0^{p,\infty }\,,\)
\(0 < p < 1\,, 0 < q \leq \infty \) and \(X\) is \(\ell^s\,,0 < s \leq \infty \). Here \(\ell^s\) is the usual Lebesgue sequence space consisting of
\(s\)-summable sequences in \(W\) when \(s \not= \infty \,,\) and bounded sequences in \(W\) when \(s = \infty \). Actually we do a little
more. If \(E\) is either \(H^{p,q}\) or \(H_0^{p,\infty }\,,\) \(0 < p < 1, 0 < q < \infty \,,\) we find \((E,X)\) whenever \(X\) is \(\ell^s\)-reflexive for
some \(0 < s \leq \infty \). We also find the multiplier spaces \((H^{p,\infty },X)\) for solid target spaces \(X\). Recall that a vector subspace
\(X\) of \(W\) is said to be solid if for every \(x = \{x_n\}, y = \{y_n\} \in W\,,\) we have \(y \in X\) whenever \(x \in X\) and \(|y_n| \leq |x_n|, n \in \mathbb{N}_0\).
Equivalently, \(X\) is solid if either of the conditions 
\begin{equation}\label{5.1}
\ell^\infty \subset(X,X)\mbox{ or }(\ell^\infty ,X)=X 
\end{equation}
\noindent
are satisfied. The notation \(s(X)=(\ell^\infty ,X)\) is commonly used. In general, \(s(X)\) is the largest solid subspace of \(X\).
We note here that for arbitrary spaces \(E\) and \(X\,,\) the multiplier space \((E,X)\) is solid whenever the target space \(X\)
is solid. Consequently \(\ell^s\)-reflexive spaces are solid. Thus the result for \(H^{p,\infty }\) is more general than the 
corresponding result for  \(H_0^{p,\infty }\).

The determination of \((E,X)\) in these cases and others frequently requires using the analytic transform to identify \((E,X)\)
with a weighted sequence space. If \(X\) is a vector subspace of \(W\) and the element \(w \in W\,,\) we define the weighted
space \(X_w = B_w^{-1}(X) = \{y \in W : wy \in X\}\). If \(X\) is a quasi-Banach space with quasinorm \(||\cdot ||_X\) then
\((X_w,||\cdot ||_{X_{w}})\)is a quasi-Banach space where \(||y||_{X_{w}} = ||yw||_X\,,\) \(y\in X_w\). For \(-\infty < \alpha < \infty \), let
\(w_\alpha =\{w_\alpha (n) : n \in \mathbb{N}_0\}\) be the power sequence defined by \(w_\alpha (0) = 1\) and \(w_\alpha = n^\alpha \) for \(n \not= 0\).
In this case we will write \((X_\alpha ,||\cdot ||_{X_\alpha })\) in place of \((X_w,||\cdot ||_{X_{w}})\). The following lemma gives the relationship between \((E,X)\)
and \((E_\alpha ,X_\beta )\). The proof is purely algebraic, [48].
\begin{Lem} Let \(E\) and \(X\) be vector subspaces of \(W\) and let \(-\infty < \alpha, \beta < \infty \). Then
\[(E_\alpha ,X_\beta ) = (E_{\alpha - \beta }, X) = (E, X_{\beta - \alpha }) = (E,X)_{\beta - \alpha }\,.\]
\end{Lem} 

Of special interest to us are the dyadically blocked sequence spaces \(\ell(p,q)\) and their weighted analogs.
These spaces are defined as follows. Let\\
\(0 < p,q \leq \infty \,,\) set \(I_0 = \{0\}\), and for \(n \in \mathbb{N}_0, n > 0\,,\) set \(I_n = \mathbb{N}_0 \cap [2^{n-1},2^n)\).
Then \(\ell(p,q)\) is the subspace of \(W\) consisting of elements \(x = \{x_n\}\) such that \(||x||_{\ell(p,q)} < \infty \,,\) where
\(||x||_{\ell(p,q)} = ||||\{x_k\}_{k\in I_n}||_{\ell^p}||_{\ell^q}\,.\)
For \(-\infty < \alpha < \infty \,,\) we write \((\ell(p,q,\alpha ), ||\cdot ||_{\ell(p,q,\alpha )})\) for the weighted space \((\ell(p,q)_\alpha ,||\cdot ||_{\ell(p,q)_\alpha})\).
The spaces \((\ell(p,q,\alpha ), ||\cdot ||_{\ell(p,q,\alpha )})\) are quasi-Banach spaces and \((\ell(p,p,0), ||\cdot ||_{\ell(p,p,0)}) = (\ell^p,||\cdot ||_{\ell^p})\), where
\(||\cdot ||_{\ell^p}\) is the standard quasinorm on \(\ell^p\). Furthermore, \(\ell(p,q,\alpha )\hookrightarrow \ell_\alpha ^\infty \hookrightarrow H(\mathbb{D})\), and
hence \(\ell(p,q,\alpha )\) is also a FK-space. The multipliers between these spaces are well-known. The following result is due mainly to Kellogg [31], see
also [26], [27], [28], and [51]. Before stating the theorem we introduce notation. Let \(0 < q, s \leq \infty \). Then \(q\ast s\)
is the extended real number defined by \(q\ast s = s\) if \(q = \infty \,,\) \(q\ast s = \frac{qs}{q -s}\) if \(0 < s < q < \infty \,,\) and
\(q\ast s = \infty \) if \(0 < q \leq s \leq \infty \).
\begin{Thm}Let \(0 < p, q, r, s \leq \infty \), \(- \infty < \alpha, \beta < \infty \). Then
\[(\ell(p,q,\alpha ),\ell(r,s,\beta )) = \ell(p\ast r, q\ast s, \beta - \alpha )\,.\]
\end{Thm}
\noindent We remark that a consequence of Theorem 5.1 is that \(\ell(p,q,\alpha )\) is \(\ell^s\)-reflexive for
\(- \infty < \alpha < \infty\) and \(0 < s \leq p, q \leq \infty \). The space \(c_0\), of null sequences, is an example
of a space which fails to be \(\ell^s\)-reflexive for every \(0 < s \leq \infty \). In addition to Theorem 5.1 we need
a lemma. Lemma 5.2 is due to Aleksandrov [2], but may also be seen to follow from Theorem 4.1. In the lemma
the spaces \(H^{p,\infty }\) and \(H_0^{p,\infty }\) are considered as sequence spaces of Taylor coefficients.
\begin{Lem} Let \(0 < p < 1\). Then\\

\noindent(i) \(H^{p,\infty }\hookrightarrow \ell_{1-1/p}^\infty \,,\)\\
(ii) \(H_0^{p,\infty }\hookrightarrow (c_0)_{1-1/p}\,.\)
\end{Lem} 
\begin{Thm}Let \(0 < p < 1\) and let \(X\) be a solid FK-space satisfying \(X\hookrightarrow H(\mathbb{D})\).Then
\((H^{p,\infty },X) = X_{1/p -1}\,.\)
\end{Thm}
\begin{proof}
Since \(X\) is solid, it follows that \(X_{1/p -1}\) is solid. Therefore, using Lemma 5.2, we find

\begin{equation}\label{5.2}
X_{1/p -1} = (\ell^\infty ,X_{1/p -1}) = (\ell_{1-1/p}^\infty,X). 
\end{equation}

\noindent
Then the inclusion

\begin{equation}\label{5.3}
X_{1-1/p} \subset (H^{p,\infty },X) 
\end{equation} 

\noindent
results from \eqref {5.2} and Lemma 5.1.\\

\noindent We obtain the reverse inclusion of \eqref{5.3} as follows. Let

\begin{equation}\label{5.4}
g(z) = (1 - z)^{-\frac{1}{p}}, z \in \mathbb{D} 
\end{equation}

\noindent
Then for \(0 < t < 1\,,\)

\begin{equation}\label{5.5}
g^*(t) \sim t^{-\frac{1}{p}}\mbox{ and } g \in H^{p,\infty }. 
\end{equation}

\noindent The Taylor coefficient sequence of \(g\) is

\begin{equation}\label{5.6}
\Big\{\frac{\Gamma (n + 1/p)}{\Gamma (1/p)n!}\Big\}. 
\end{equation}

\noindent A well-known consequence of Stirling's formula is that 
\begin{equation}\label{5.7}
\Big\{\frac{n^{1/p -1}\\ n!}{\Gamma (n + 1/p)}\Big\}\in \ell^\infty .
\end{equation}
\noindent
Therefore we deduce the reverse inclusion of \eqref{5.3} from \eqref{5.1} and \eqref{5.5} through \eqref{5.7}.
\end{proof}
\begin{Cor} Let \(0 < p < 1, 0 < s \leq \infty \). Then \((H^{p,\infty },\ell^s)=\ell_{1/p - 1}^s\).
\end{Cor}

\begin{Thm}Let \(0 < p < 1, 0 < s \leq \infty \). Then \((H_0^{p,\infty },\ell^s)=\ell_{1/p - 1}^s\).
\end{Thm}
\begin{proof}
We prove only the case \(s\not= \infty \), the other case being similar.
Using Corollary 5.1 and Lemma 3.1(i) we obtain

\begin{equation}\label{5.8}
\ell_{1/p - 1}^s = (H^{p,\infty }, \ell^s) \subset (H_0^{p,\infty }, \ell^s).
\end{equation}

Next we show the reverse inclusion of \eqref {5.8} holds. Fix \(\lambda = \{\lambda _n\}\in(H_0^{p,\infty },\ell^s)\). Since
\(\ell^s\) is solid there is no loss of generality in assuming \(\lambda _n \geq 0\), \(n \in \mathbb{N}_0\). Furthermore, the solidity of \(\ell_{1/p -1}^s\) and \eqref{5.7},
show that \(\lambda \in \ell_{1/p-1}^s\) if and only if \(\{\frac{\Gamma (n+1/p)\lambda _n}{\Gamma(1/p)n!}\}\in \ell^s\).
Consider the operator \(B_\lambda  \in \cal L\)\((H_0^{p,\infty }, \ell^s)\) corresponding to \(\lambda \). That is,
for \(f \in H_0^{p,\infty }\,,\) with Taylor series representation \eqref {2.3} we have 

\begin{equation}\label{5.9}
B_\lambda (f) = a\lambda = \{a_n \lambda _n\}\,, 
\end{equation}  

and

\begin{equation}\label{5.10}
||B_\lambda ||_{ {\cal L}(H_0^{p,\infty },\ell^s) } = ||\lambda ||_{(H_0^{p,\infty },\ell^s)} < \infty . 
\end{equation}

\noindent Let \(g\) be the Cauchy-type function in \eqref{5.4}. Then for \(0 < r, t < 1\,,\)
 
\begin{equation}\label{5.11}
(g_r)^*(t)\sim C_1(1 -r)^{-\frac{1}{p}}\chi _{[0,1- r)}(t) + C_2t^{-\frac{1}{p}}\chi _{(1- r,1)}(t) 
\end{equation}

\noindent where the constants \(C_j\), \(j = 1, 2\) are independent of r and t, and \(\chi _A\) denotes the 
characteristic function of a set \(A\). Thus \eqref{5.5} and \eqref {5.11} imply

\begin{equation}\label{5.12}
\sup_{0 \leq r < 1}||g_r||_{H^{p,\infty }}\leq C||g||_{H^{p,\infty }}< \infty . 
\end{equation} 
For \(0 \leq r < 1\), put \(\Phi (r) = ||B_\lambda (g_{r^{1/s}})||_{\ell^s}^s\). It follows from \eqref {5.10} and \eqref {5.12} that \(\Phi \)
is bounded on \([0,1)\).
From \eqref {5.6} and \eqref {5.9} we see that \(g_{r^{1/s}} = \{\frac{\Gamma (n+1/p)\lambda _n r^{n/s}}{\Gamma (1/p)n!}\}\) 
and \(\Phi (r) = ||\{\frac{\Gamma (n + 1/p) \lambda _n r^{n/s}}{n!\Gamma (1/p)}\}||_{\ell^s}^s\).  
It therefore follows that \(\Phi \) is an increasing function of \(r\).  
From these observations we deduce that \(\lim_{r\rightarrow 1^-}\Phi (r)\) exists, hence the positive sequence 
\(\{(\frac{\Gamma (n+1/p)\lambda _n}{n! \Gamma (1/p)})^s\}\)is Abel summable and consequently
belongs to \(\ell^1\). But then \(\{\frac{\Gamma (n+1/p)\lambda _n}{n!\Gamma (1/p)}\} \in \ell^s\) which is what we needed to show.   
\end{proof}

The space \(\ell^s\) is a rearrangement invariant quasi-Banach function space with \(\mathbb{N}_0\) equipped with counting measure
as the underlying measure space. Quasi-Banach function spaces are always solid. A rearrangement invariant quasi-Banach function space \(X\) is called maximal
if every quasinorm bounded increasing sequence in \(X\) is bounded above in \(X\), see [30]. It is not hard to see that if \(X\) is a maximal rearrangement 
invariant quasi-Banach function space then \((H_0^{p,\infty },X) = X_{1/p - 1}\). Another situation where we have    
\((H_0^{p,\infty },X) = X_{1/p - 1}\) is when \(X\) is \(\ell^s\)-reflexive. For \(0 < s \leq \infty \) and an arbitrary space \(X\) we use the notation \(X^{K(s)}\) for \((X,\ell^s)\) and \(X^{K(s)K(s)}\) for \((X^{K(s)})^{K(s)}\).
Thus the space \(X\) is \(\ell^s\)-reflexive if and only if \(X^{K(s)K(s)} = X\). We have the following generalization of Theorem 5.3.
\begin{Thm}
Let \(0 < p < 1\) and let \(X\) be a FK-space which is \(\ell^s\)-reflexive for some \(0 < s \leq \infty \). Then\\

\noindent
\((H_0^{p,\infty },X) = X_{1/p -1}\,.\)   
\end{Thm}
\begin{proof}
Since \(X\) is \(\ell^s\)-reflexive, \(X_{1/p-1}\) is solid. Hence Theorem 5.2 and Lemma 3.1(i) produce the inclusion
\begin{equation}\label{5.15}
X_{1/p-1} = (H^{p,\infty },X)\subset (H_0^{p,\infty },X).  
\end{equation}

To get the reverse inclusion of \eqref{5.15} we use Lemma 3.1(ii), the \(\ell^s\)-reflexivity of \(X\,,\) and the identity

\begin{equation}\label{5.16}
(\ell_{1/p -1}^s)^{K(s)} = \ell_{1 - 1/p}^\infty\,. 
\end{equation}

\noindent
Then using \eqref{5.2}, \eqref{5.15}, Theorem 5.3, Lemma 3.1(ii), Lemma 5.1, and \eqref{5.16} we get

\begin{eqnarray*}
(H_0^{p,\infty },X) & \subset & (X^{K(s)},(H_0^{p,\infty })^{K(s)})\\
& = & (X^{K(s)},\ell_{1/p-1}^s)\\
& \subset & (\ell_{1 -1/p}^\infty , X^{K(s)K(s)})\\
& = & (\ell_{1-1/p}^\infty ,X)\\
& = & X_{1/p-1}\\
& \subset & (H_0^{p,\infty },X)\,.
\end{eqnarray*}
and the proof is complete.
\end{proof}

\noindent
In Theorem 5.4, the hypothesis that \(X\) be \(\ell^s\)-reflexive for some \(0 < s \leq \infty \) cannot be omitted.
\begin{Cor}
Let \(0 < p < 1\). Then \((H_0^{p,\infty },c_0) = \ell_{1/p -1}^\infty \,.\)
\end{Cor}

\noindent
Proof: Use Lemma 3.1(i) and (ii), Lemma 5.2(ii), Theorem 5.4 and the identity \(\ell_{1/p -1}^\infty = ((c_0)_{1-1/p},c_0)\). \\ 

We turn now to the study of the multiplier space \((H^{p,q},X)\) where \(X\) is a \(\ell^s\)-reflexive FK-space.
We need two results from the theory of mixed norm spaces. The first of these is part of the folklore. The case \(q = t\) may be found in [8].
\begin{Lem}
Let \(0 < p \leq 2, 0 < q \leq t \leq \infty , 0 < \alpha  < \infty , - \infty < \beta < \infty \). Set \(p_0 = min(1,p)\) and \(p_1 = max(1,p)\). 
Let \(p_1'\) be the H\"{o}lder conjugate of \(p_1\), \(1/p_1 + 1/p_1' =1\). Then there is the embedding
\[H(p, q, \alpha , \beta ) \hookrightarrow \ell(p_1', t, 1 - 1/p_0 + \beta - \alpha )\,.\] 
\end{Lem}

\noindent
The second result we need is a theorem of Pavlovic characterizing the multiplier spaces \((H(p, q, \alpha ), \ell^s)\)
for \(0 < p \leq 1, 0 < q, s \leq \infty , 0 < \alpha < \infty \,,\) [43]. See also [28]. For some special cases of the theorem see [1]
and [39].
\begin{Thm}
Let \(0 < p \leq 1, 0 < q, s \leq \infty , 0 < \alpha < \infty \). Then
\[(H(p, q, \alpha ), \ell^s) = \ell(s, q\ast s, \alpha + 1/p - 1)\,.\] 
\end{Thm}

\noindent
Lemma 2.3 may be used to obtain the following extension of Theorem 5.5.
\begin{Cor}
Let \(0 < p \leq 1, 0 < q, s \leq \infty , 0 < \alpha < \infty , - \infty < \beta < \infty\). Then
\[(H(p, q, \alpha, \beta  ), \ell^s) = \ell(s, q\ast s, \alpha - \beta + 1/p - 1)\,.\] 
\end{Cor}
\noindent
Duren and Shields showed that \((H^p, \ell^s) = \ell(s, \infty , 1/p -1)\) for \(0 < p < 1, \\ p \leq s \leq \infty\), [10], [12], [13].
In [27] Jevtic and Pavlovic showed that\\
\((H^p, \ell^s) = \ell(s,p\ast s, 1/p -1)\) for the case \(0 < s < p < 1\). Theorem 5.6
below extends these results to the Hardy-Lorentz space setting.  
\begin{Thm}
Let \(0 < p < 1, 0 < q < \infty , 0 < s \leq \infty \). Then

\begin{equation}\label{5.17}
(H^{p,q},\ell^s) = \ell(s,q\ast s,1/p -1)\,.
\end{equation}
\end{Thm}
\begin{proof}
From Theorem 4.1 we have the embeddings

\begin{equation}\label{5.18}
H(p_0,q,\beta + 1/p - 1/p_0,\beta ) \hookrightarrow H^{p,q} \hookrightarrow H(1,q,1/p -1)\,. 
\end{equation}

\noindent
where \(\beta > 1/p_0 - 1/p > 0\). Applying Lemma 3.1(i) to \eqref{5.18} and using Corollary 5.3
we have 

\begin{eqnarray*}
\ell(s,q\ast s,1/p -1) & = & (H(1,q,1/p -1),\ell^s)\\
& \subset & (H^{p,q},\ell^s)\\
& \subset & (H(p_0,q,,\beta + 1/p -1/p_,\beta )\\
& = & \ell(s,q\ast s,1/p -1)\,,\\ 
\end{eqnarray*}

\noindent
which establishes \eqref{5.17}.
\end{proof}

\noindent Theorem 5.6 generalizes to
\begin{Thm}
Let \(0 < p < 1\), \(0 < q < \infty \). Let \(X\) be a FK-space which is \(\ell^s\)-reflexive for some \(0 < s \leq \infty \)
and set  \(t = max(q,s)\). Then 
\[(H^{p,q},X) = (\ell(\infty ,t,1-1/p),X)\,.\] 
\end{Thm}
\begin{proof}
Observe that for \(0 < q, s \leq \infty \),

\begin{equation}\label{5.19}
(q\ast s)\ast s = t\,. 
\end{equation}

\noindent
Now apply Lemma 3.1(iii) followed by Theoren 5.6 to obtain

\begin{eqnarray}\label{5.20}
(H^{p,q},X) & \subset & (X^{K(s)},(H^{p,q})^{K(s)})\notag\\ 
& = & (X^{K(s)},\ell(s,q\ast s,1/p -1))\,.\notag\\
\end{eqnarray}

\noindent
Since \(X\) is \(\ell^s\)-reflexive, a second application of Lemma 3.1(iii) to the last space in \eqref{5.20}
together with Theorem 5.1 and \eqref{5.19} yields

\begin{eqnarray}\label{5.21}
(X^{K(s)},\ell(s,q\ast s,1/p -1) & \subset & (\ell(s,q\ast s,1/p -1)^{K(s)},X^{K(s)K(s)})\notag\\
& = & (\ell(\infty , (q\ast s)\ast s,1-1/p),X)\notag\\
& = & (\ell(\infty ,t,1-1/p),X)\,.\notag\\
\end{eqnarray}

Then, starting with the last space in \eqref{5.21}, use Lemma 3.1(iii) three times, first with Lemma 5.3, then with Theorem 4.1,
and finally with the Hardy-Lorentz analog of inclusion \eqref{2.1}. As a result we get

\begin{eqnarray}\label{5.22}
\ell(\infty ,t,1 - 1/p),X) & \subset & (H(1,t,1/p - 1),X)\notag\\
& \subset & (H^{p,t},X)\notag\\
& \subset & (H^{p,q},X)\,.\notag\\
\end{eqnarray}
Combining \eqref{5.20} through \eqref{5.22} completes the proof.
\end{proof}

Theorem 5.4 and Theorem 5.7 may be used to compute the multiplier spaces \((E, ces(s))\), where \(E\) is one of the Hardy-Lorentz spaces 
\(H_0^{p,\infty }\) or \(H^{p,q}\), \(0 < p < 1\), \(0 < q \leq \infty \) and for \(1 < s < \infty \), \(ces(s)\) is the Cesaro sequence space consisting of 
sequences \(\{x_k\} \in W\) satisfying \(\sum_{n = 1}^{\infty } (\frac{1}{n} \sum_{k=1}^{n} |x_k|)^s < \infty  \). Since it is known that \(ces(s) = \ell(1,s,\frac{1}{s} - 1)\),
[20], then we have the following result.

\begin{Cor}
Let \(0 < p < 1\), \(0 < q \leq \infty \), \(1 < s < \infty \). Then\\
\noindent (i) \((H^{p,q},ces(s)) = \ell(1,q\ast s, 1/p + 1/s - 2)\),\\
(ii) \(H_0^{p,\infty },ces(s)) = \ell(1, s, 1/p + 1/s - 2)\).
\end{Cor}
\begin{Cor}
Let \(0 < p < 1\) and suppose \(X\) is a FK-space which is \(\ell^s\)-reflexive for some \(0 < s \leq \infty \). Then for every
\(0 < q \leq s\), 
\[(H^{p,q},X) = (H^{p,s},X) = (\ell(\infty ,s,1 - 1/p),X)\,.\]
\end{Cor}
\noindent
Corollary 5.5 asserts that given a FK-space \(X\) which is \(\ell^s\)-reflexive for some \(0 < s \leq \infty \) and a number \(0 < p < 1\,,\) the \(X\)-valued
multiplier spaces for the Hardy-Lorentz space scale \(\{H^{p,q}\}_{0 < q \leq s}\) will coincide. This is really due to the fact that the Bergman-Sobolev spaces
appearing in the proof of Theorem 5.5 enjoy this property. We conclude this section with a Hardy-Lorentz analog of a well-known result of Hardy and
Littlewood for \(H^p\) spaces.
\begin{Cor}
Let \(0 < p < 1\), \(0 < q < \infty \). Suppose \(f \in H^{p,q}\) has Taylor series representation \(f(z) = \sum_{n = 0}^{\infty } a_n z^n\,,z \in \mathbb{D}\).
Then 

\begin{equation}\label{5.23}
\Big(\sum_{n = 1}^{\infty } n^{q(1 -\frac{1}{p}) - 1}|a_n|^q\Big)^{1/q}\leq C ||f||_{H^{p,q}}\,.  
\end{equation}

\end{Cor} 

\noindent
The diagonal case \(p = q\) is a result of Hardy and Littlewood and is actually valid for \(0 < p \leq 2\). 
We also mention that for \(f\in H^p\), \(0 < p \leq 1, p < q < \infty\,,\) the series \eqref {5.23} is known to converge, [12].
Since \(\{n^{1 - \frac{1}{p} -\frac{1}{q}}\}\) belongs to \(\ell(q,\infty ,1/p - 1)\,,\) Corollary 5.5 follows from the identification \((H^{p,q},\ell^q) = \ell(q,\infty ,1/p - 1)\).
Finally we note that \eqref {5.23} remains valid if \(H^{p,q}\) is replaced by \(H(1,q,1/p -1)\).

\setcounter{equation}{0}
\section{Multipliers of \(H^{p,q}\) and \(H_0^{p,\infty }\,,\) \(0 < p < 1, \\0 < q < \infty \) into Bergman-Sobolev spaces}
Our main tools in this section are Theorem 4.1 and the following result of Pavlovic [43], see also [28].
\begin{Thm}
Let \(0 < q, s, t \leq \infty , 0 < p \leq min(1,s), 0 < \alpha, \beta < \infty \). Then

\begin{equation}\label{6.1}
(H(p, q, \alpha ), H(s, t, \beta )) = H(s, q\ast t, \beta , \alpha + 1/p -1)\,. 
\end{equation}

\end{Thm}
\noindent
Note that by Lemma 2.2, the space on the right-hand side of \eqref{6.1} coincides with \(H(s,q\ast t,1,1/p + \alpha - \beta )\).
\begin{Cor}
Let \(0 < q, s, t \leq \infty , 0 < p \leq min(1,s), 0 < \alpha, \beta < \infty ,\\ - \infty < \delta , \gamma < \infty \). Then\\

\noindent  
(i) \((H(p,q,\alpha, \delta ),H(s,t,\beta ,\gamma )) = H(s,q\ast t,1,1/p + \alpha - \beta + \gamma - \delta )\,,\)\\ 
(ii) \((H_0(p,\infty ,\alpha, \delta),H(s,t,\beta, \gamma)) = H(s,t,1,1/p + \alpha -  \beta + \gamma - \delta)\,,\)\\
(iii) \((H_0(p,\infty ,\alpha, \delta),H_0(s,\infty ,\beta, \gamma)) = H(s,\infty ,1,1/p + \alpha -  \beta + \gamma - \delta)\,,\)\\
(iv)\((H(p,\infty ,\alpha, \delta),H_0(s,\infty ,\beta, \gamma)) = H_0(s,\infty ,1,1/p + \alpha -  \beta + \gamma - \delta)\,,\)\\
(v) \((H(p,q,\alpha, \delta),H_0(s,\infty ,\beta, \gamma)) = H(s,\infty ,1,1/p + \alpha -  \beta + \gamma - \delta)\) if \(q \not= \infty \,.\)\\
\end{Cor}
\begin{proof}
(i) Let \(g\in (H(p,q,\alpha, \delta ),H(s,t,\beta ,\gamma ))\). By Lemma 2.3 the maps \(D_\delta \) and \(D^\gamma \) are
continuous surjective isomorphisms
\[D_\delta : H(p,q,\alpha )\rightarrow H(p,q,\alpha, \delta)\]
\[D^\gamma : H(s,t,\beta, \gamma)\rightarrow H(s,t,\beta )\]
and hence \(D^\gamma \circ g \circ D_\delta\in (H(p,q,\alpha ),H(s,t,\beta ))\). But \(D^\gamma \circ g \circ D_\delta \) is the
same multiplier as \((g^{[\gamma ]})_{[\delta ]}\). Therefore, from \eqref{6.1}, Lemmas 2.2 and 2.3, we deduce

\begin{eqnarray}\label{6.2}
(H(p,q,\alpha, \delta ),H(s,t,\beta ,\gamma )) & \subset & D_\gamma (D^\delta ((H(p,q,\alpha ),H(s,t,\beta))))\notag\\
& = & D_\gamma (D^\delta (H(s,q\ast t,\beta ,\alpha +1/p -1)))\notag\\
& = & D_\gamma (H(s,q\ast t,\beta ,\alpha - \delta + 1/p -1))\notag\\
& = & H(s,q\ast t,\beta,\alpha + \gamma - \delta + 1/p -1)\,.\notag\\ 
\end{eqnarray}

\noindent
Similarly,

\begin{equation}\label{6.3}
(H(p,q,\alpha ),H(s,t,\beta ))\subset D_\delta (D^\gamma (H(p,q,\alpha, \delta), H(s,t,\beta, \gamma))\,.   
\end{equation}

\noindent
Therefore \eqref{6.3}implies

\begin{eqnarray}\label{6.4}
H(s,q\ast t,\alpha + \gamma - \delta + 1/p -1) & = & D_\gamma (D^\delta (H(s,q\ast t,\beta ,\alpha + 1/p -1)))\notag\\
& = & D_\gamma (D^\delta ((H(p,q,\alpha ),H(s,t,\beta ))))\notag\\ 
& \subset & (H(p,q,\alpha, \delta),H(s,t,\beta ,\gamma ))\,.\notag\\ 
\end{eqnarray}

\noindent
Then \eqref{6.2} and \eqref {6.4} together give us (i).\\
(ii) This follows from (i) and the monotonicity of the means \(M_k(\cdot ,f)\) for \(f\in H(\mathbb{D}), 0 < k \leq \infty \).\\
(iii) This follows from (ii) and the fact that for every \(0 < k \leq \infty ,\\ 0< \eta < \infty , - \infty < \nu < \infty\,,\) the function \(F\in H(k,\infty ,\eta, \nu)\)
belongs to \( H_0(k,\infty ,\eta, \nu)\) if and only if \(F_r\rightarrow F\) in \(H(k,\infty ,\eta, \nu)\) as \(r\rightarrow1^-\).\\
(iv) First note that the inclusion

\begin{equation}\label{6.5}
(H_0(s,\infty ,1,1/p + \alpha - \beta + \gamma - \delta) \subset (H_0(p,\infty ,\alpha,\gamma),H(s,\infty ,\beta, \gamma)) 
\end{equation}

\noindent
follows from (i) and the fact that for every \(0 < k \leq \infty\,, 0< \eta < \infty\,,\,\\ - \infty < \nu < \infty\,,\) the function \(F\in H(k,\infty ,\eta, \nu)\)
belongs to \( H_0(k,\infty ,\eta, \nu)\) if and only if \(F_r\rightarrow F\) in \(H(k,\infty ,\eta, \nu)\) as \(r\rightarrow1^-\).\\
\indent
To prove the reverse inclusion of \eqref{6.5} we consider the case \(\delta = \gamma = 0\) first. Let \(g\in (H(p,\infty ,\alpha ),H_0(s,\infty ,\beta ))\).
The Cauchy-type function \(F(z) = (1 - z)^{-\alpha - \frac{1}{p}}, z \in \mathbb{D}\) belongs to \(H(p,\infty, \alpha)\). 
Therefore \(g\ast F\in H_0(s,\infty ,\beta )\). But for any \(w\in \mathbb{D}\), \(g \ast F(w) = \Gamma (\alpha + 1/p)^{-1}g^{[\alpha + \frac{1}{p} -1]}(w)\).
Hence\\ \(g\in H_0(s,\infty ,\beta ,\alpha + 1/p -1)\) which is equivalent to \(g\in H_0(s,\infty ,1,1/p + \alpha -  \beta)\) by Lemma 2.2.
For the general case we argue as in the proof of (i), using Lemma 2.3 to write
\[(H(p,\infty ,\alpha, \delta), H_0(s,\infty ,\beta, \gamma)) = D^\delta (D_\gamma ((H(p,\infty, \alpha),H_0(s,\infty, \beta))).\]
Then use the validity of (iv) for the case \(\delta = \gamma = 0\).\\
(v) The proof is similar to the proof of (iv).
\end{proof}
\begin{Thm}
Let \(0 < q, s, t \leq \infty , 0 < \beta < \infty , 0 < p < min(1,s),\\ - \infty < \gamma < \infty \). Then\\

\noindent
(i) \((H^{p,q},H(s,t, \beta, \gamma )) = H(s,q\ast t, 1, 1/p + \gamma - \beta),\)\\ 
(ii) \((H^{p,q},H_0(s,\infty ,\beta, \gamma)) = H(s,\infty ,1,1/p + \gamma - \beta)\,, q \not= \infty \)\\  
(iii) \((H^{p,\infty },H_0(s,\infty ,\beta, \gamma)) = H_0(s,\infty ,1,1/p + \gamma - \beta),\)\\
(iv) \((H_0^{p,\infty },H(s,t,\beta, \gamma)) = H(s,t,1,1/p + \gamma - \beta),\)\\
(v) \((H_0^{p,\infty },H_0(s,\infty ,\beta, \gamma)) = H(s,\infty ,1,1/p + \gamma - \beta).\)
\end{Thm} 
\noindent
\begin{proof}
We prove (i) only, the proofs of (ii) through (v) being similar. Let \(0 < p_0 < p\) and \(\delta > 1/p_0 - 1/p\). 
By Theorem 4.1 we have embeddings

\begin{equation}\label{6.6}
H(p_0,q,\delta + 1/p - 1/p_0, \delta )\hookrightarrow H^{p,q} \hookrightarrow H(s,q, 1/p - 1/s), 
\end{equation}

\noindent     
By Corollary 6.1, the multiplier spaces \((H(p_0,q,\delta + 1/p - 1/p_0, \delta ), H(s,t,\beta, \gamma))\) and \((H(s,q, 1/p - 1/s),H(s,t,\beta, \gamma))\)
are both equal to\\
\(H(s,q\ast t, 1, 1/p + \gamma - \beta)\). Then (i) follows from this fact and \eqref{6.6}.
\end{proof}

\setcounter{equation}{0}
\section{Multipliers of \(H^{p,q}\) and \(H_0^{p,\infty }\,,\) \(0 < p < 1,\\ 0 < q \leq \infty \) into Hardy spaces} 
In this section we consider the multiplier spaces \((H^{p,q},H^s)\) and \((H_0^{p,\infty },H^s)\) for \(0 < p < 1\), \(0 < q \leq \infty\), \(0 < s < \infty \).
We are able to determine these spaces for the cases \(0 < q \leq \mbox{ min}(2,s)\) and \(0 < q \leq \infty , s = 2\). 
Since \(H^2 = \ell^2\,,\) the second case was addressed in Section 5. We restate that result in terms of Bergman-Sobolev spaces below. To
do this we need the following lemma from [8], see also [40], [43], [53].
\begin{Lem}
Let \(0 < q \leq \infty , 0 < \alpha < \infty , - \infty < \beta < \infty .\) Then
\[H(2,q,\alpha ,\beta ) = \ell(2,q,\beta - \alpha ).\]
\end{Lem}

\noindent
Using Lemma 7.1, the identification \(H^2 = \ell^2\), and either Corollary 5.1 or Theorem 5.6, we have the following result.
\begin{Thm}
Let \(0 < p < 1, 0 < q \leq \infty \). Then
\[(H^{p,q},H^2) = H(2,q\ast 2,1,1/p).\]
\end{Thm}

\noindent
We turn now to the case \(0 < q \leq\) min\((2,s)\). First we record the following.
\begin{Thm}
Let \(0 < s \leq \infty , 0 < p < min(1,s)\). Then
\[(H^p,H^s) = H(s,\infty ,1,1/p).\]
\end{Thm}

\noindent
Theorem 7.2 dates back to Hardy and Littlewood who observed that\\ \(H(s,\infty ,1,p) \subset (H^p,H^s)\) for 
\(0 < p < 1 \leq s \leq \infty\), [20], [21]. The case \(s = \infty \) corresponding to the Duren-Romberg-Shields
Theorem of [11] reduces to \eqref{4.3}. Duren and Shields proved Theorem 7.2 for the case\\ \(0 < p < 1 \leq s < \infty\), 
[13]. The proof for the case \(0 < p < s \leq 1\) is due to Mateljevic and Pavlovic, [38]. Theorem 7.3 below
represents an extension of Theorem 7.2 to \(H^{p,q}\) for \(0 < q \leq \) min\((2,s)\). We will need the following lemma.
Before stating this result we introduce some notation. For \(0 < s < \infty \), the Dirichlet-type space \({\cal {D}}^s\) is
defined to be the space \(H(s,s,1,1)\).
\begin{Lem}
Let \(0 < s \leq 2 \leq t < \infty \). Then \\

\noindent
(i) \({\cal{D}}^s \hookrightarrow H^s \hookrightarrow H(s,2,1,1)\,,\)\\
(ii) \(H(t,2,1,1)\hookrightarrow H^t\hookrightarrow {\cal{D}}^t\,.\)     
\end{Lem}
\noindent
For statements and proofs of Lemma 7.2 the reader may consult [8],[17],[28],[34],[36], and [41].
Recently, A. Baernstein, D. Girela, and J. A. Pelaez have shown that for all \(0 < s < \infty \), \(H^s\cap\ {\cal {U}} = {\cal{D}}^s\cap\cal U\),
where \(\cal U\) is the class of univalent functions on \(\mathbb{D}\), [4].
\begin{Thm}
Let \(0 < s < \infty\), \(0 < p < min(1,s),\mbox{and } 0 < q \leq min(2,s)\). Then
\[(H^{p,q}, H^s) = H(s,\infty ,1,1/p)\,.\]
\end{Thm}
\begin{proof}
Proof: Assume first that \(0 < s \leq 2\). Then using Lemma 7.2 and Lemma 3.1(ii) we have
 
\begin{equation}\label{7.1}
(H^{p,q},{\cal{D}}^s)\subset (H^{p,q},H^s)\subset (H^{p,q},H(s,2,1,1)) . 
\end{equation}

\noindent
By Theorem 6.2(i), both of the endpoint spaces in \eqref{7.1} are equal to \(H(s,\infty ,1,1/p)\).
\noindent
For the case \(2 \leq s < \infty \), the reverse inclusion of \eqref{7.1} holds and the rest of the proof is exactly the same as in the first case.
\end{proof}

\noindent
Theorem 7.3 has the following corollary.
\begin{Cor}
Let \(0 < q \leq s \leq 2, 0 < p < min(1,s)\). Then

\begin{equation}\label{7.2}
(H^{p,q},H^{s,q}) = \bigcap _{q \leq t \leq s}(H^{p,t},H^{s,t})
\end{equation}
\end{Cor}

\begin{proof}
Proof: We prove the inclusion \((H^{p,q},H^{s,q})\subset \bigcap_ {\,q \leq t \leq s} (H^{p,t},H^{s,t})\) with
the reverse conclusion being obvious. Let \(g\in (H^{p,q},H^{s,q})\). Since \(q \leq s \leq 2\,,\) we have
\((H^{p,q},H^{s,q})\subset (H^{p,q},H^s) = (H^{p,s},H^s)\) by the Hardy-Lorentz analog of \eqref{2.1},
Lemma 3.1(i) and Theorem 7.3. Thus g is a bounded multiplier for 

\[g : H^{p,s}\rightarrow H^s\mbox{ and } g : H^{p,q}\rightarrow H^{s,q}\,.\] 

\noindent
Therefore, by interpolation, we find g is also bounded as a multiplier

\[g : (H^{p,q},H^{p,s})_{\theta ,t}\rightarrow (H^{s,q},H^s)_{\theta ,t}\,,\]

\noindent
for \(0 < \theta < 1\,,\) and \(\frac{1}{t} = \frac{1 - \theta }{q} + \frac{\theta }{s}\). 
Then an application of Theorem 2.1(ii) implies \(g\) is bounded as a multiplier 

\[g : H^{p,t}\rightarrow H^{s,t}\mbox{ for all } q \leq t \leq s\,.\] 
\end{proof}
\setcounter{equation}{0}
\section{Multipliers of \(H^{p,q}\) and \(H_0^{p,\infty }, 0 < p < 1, \\0 < q \leq \infty \) into analytic
Lipschitz spaces, analytic Zygmund spaces, Bloch spaces, and \(BMOA\)}
In this section we apply the results of the previous two sections to some specific target spaces
belonging to the class of Bergman-Sobolev spaces. The target spaces we have in mind are the analytic 
Lipschitz and Zygmund spaces and the Bloch spaces. We also have some results for the case when the 
target space is \(BMOA\). For the discussion that follows we assume \(f\in H(\mathbb{D})\) and that
\(f\) has non-tangential limits \(m\)-a.e. on \(\mathbb{T}\). We denote the resulting boundary value function by the 
same symbol \(f\). For \(1 \leq s \leq \infty \,,\) the moduli of continuity \(\omega _s (f)(t)\) and \(\Omega _s (f)(t)\) of
\(f\) are defined for \(t > 0\) by
\(\omega _s (f)(t) = \sup_{0 < |h| \leq t}||T_h(f) - f||_s\) and \(\Omega _s (f)(t) = \sup_{0 < |h| \leq t}||T_h(f) - 2f + T_{-h}(f)||_s\),
where \(T_h\) is the translation operator given by \(T_h(f)(e^{i\theta }) = f(e^{i(\theta + h)})\).
Let \(0 < \alpha \leq 1\). Then \(f\) is said to belong to the analytic Lipschitz space \(\Lambda _\alpha ^s(\mathbb{D})\)
(resp. \(\lambda _\alpha ^s(\mathbb{D})\) if \(\omega _s (f)(t) = O(t^\alpha )\) (resp. \( o(t^\alpha )\)) as \(t\rightarrow 0^+\).
If the boundary value function \(f\in C(\mathbb{T})\) and \(\Omega _\infty (f)(t) = O(t)\) (resp. \(o(t)\)) as \(t\rightarrow 0^+\,,\)
then \(f\) is said to belong to the analytic Zygmund space \(\Lambda _\ast ^\infty (\mathbb{D})\) (resp. 
\(\lambda _\ast ^\infty (\mathbb{D})\)). For \(1 \leq s < \infty \,,\) \(f\) is said to belong to the analytic Zygmund space
\(\Lambda _\ast ^s (\mathbb{D})\) (resp.  \(\lambda _\ast ^s(\mathbb{D})\)) if \(\Omega _s (f)(t) = O(t^\alpha )\) (resp.
\(o(t^\alpha ))\) as \(t\rightarrow 0^+\). Theorem 8.1 below is a collection of well-known results of Hardy and Littlewood
identifying various analytic Lipschitz and Zygmund spaces as Bergman-Sobolev spaces. See [10] and [57]. In order to 
cover the case \(\alpha = 1\), we recall that for \(0 < s \leq \infty \,,\) \(0 < \beta < \infty \,,\) the Hardy-Sobolev space
\(H_\beta ^s = \{f \in H(\mathbb{D}) : f^{[\beta ]} \in H^s \}\).

\begin{Thm}
Let \(0 < \alpha < 1 \leq s \leq \infty \). Then\\

\noindent
(i) \(\Lambda _\alpha ^s(\mathbb{D}) = H(s,\infty ,1 - \alpha , 1)\) and \(\lambda _\alpha ^s(\mathbb{D}) = H_0(s,\infty ,1 - \alpha , 1)\,,\)\\
(ii) \(\Lambda _\ast ^s(\mathbb{D}) = H(s,\infty ,1 - \alpha , 2)\) and \(\lambda _\ast ^s(\mathbb{D}) = H_0(s,\infty ,1 - \alpha , 2)\,,\)\\
(iii) \(\Lambda _1^s(\mathbb{D}) = H_1^s\,.\)
\end{Thm} 

\noindent
We combine Theorem 8.1 with the Duren-Romberg-Shields Theorem to find the multipliers from \(H^{p,q}\) into the analytic
Lipschitz spaces \(\Lambda _\alpha ^\infty (\mathbb{D})\,,\) \(\lambda _\alpha ^\infty (\mathbb{D})\) and analytic
Zygmund spaces  \(\Lambda _\ast ^\infty (\mathbb{D})\,,\) \(\lambda _\ast ^\infty (\mathbb{D})\,,\) \(0 < \alpha , p < 1\,,\) \(0 < q \leq \infty \).

\begin{Cor}
Let \(0 < \alpha , p < 1\,,\) \(0 < q \leq \infty \). Then\\

\noindent
(i) \((H^{p,q},\Lambda _\alpha ^\infty (\mathbb{D})) = (H_0^{p,\infty },\Lambda _\alpha ^\infty (\mathbb{D})) = H(\infty ,\infty ,1, \frac{1}{p} + \alpha ) = (H^{\frac{p}{1 + \alpha p}})^\ast \,,\)\\
(ii) For \(q \not= \infty \,,\)\\ 
\[(H^{p,q},\lambda _\alpha ^\infty (\mathbb{D})) =  (H_0^{p,\infty },\lambda _\alpha ^\infty (\mathbb{D})) = H(\infty ,\infty ,1, 1/p + \alpha )=(H^{\frac{p}{1 + \alpha p}})^\ast \,,\]\\
\noindent 
(iii) \((H^{p,\infty },\lambda _\alpha ^\infty (\mathbb{D})) = H_0(\infty ,\infty ,1,\frac{1}{p} + \alpha )\,,\)\\
(iv) \((H^{p,q},\Lambda _\ast ^\infty (\mathbb{D})) = (H_0^{p,\infty },\Lambda _\ast ^\infty (\mathbb{D})) = H(\infty ,\infty ,1, \frac{1}{p} + 1) = (H^{\frac{p}{1 + p}})^\ast \,,\)\\
(v) For \(q \not= \infty \,,\)\\ 
\[(H^{p,q},\lambda _\ast ^\infty (\mathbb{D})) =  (H_0^{p,\infty },\lambda _\ast ^\infty (\mathbb{D})) = H(\infty ,\infty ,1, 1/p + 1)=(H^{\frac{p}{1 + p}})^\ast \,,\]\\
(vi) \((H^{p,\infty },\lambda _\ast ^\infty (\mathbb{D})) = H_0(\infty ,\infty ,1,\frac{1}{p} + 1)\,.\)\\
\end{Cor} 

\noindent
Corollary 8.1 shows that for \(0 < \alpha < 1\,,\) the secondary index \(q\) is irrelevant with respect to the multiplier
spaces \((H^{p,q},E)\) for the target spaces\\
\(E = \Lambda _\alpha ^\infty(\mathbb{D})\) or \(E = \Lambda _\ast ^\infty(\mathbb{D})\). The same is true for the target spaces
\(E = \lambda _\alpha ^\infty(\mathbb{D})\) or \(E = \lambda _\ast ^\infty(\mathbb{D})\) provided \(q \not= \infty \). A similar
phenomenon occurs when the target spaces are the Bloch spaces. Let us recall that the Bloch space \(\cal B\) and the little
Bloch space \({\cal B}_0\) are realized as Bergman-Sobolev spaces using the identifications \({\cal B} = H(\infty, \infty,1,1)\) and 
\({\cal B}_0 = H_0(\infty, \infty,1,1)\). The Bloch space analog of Corollary 8.1 is the following.

\begin{Cor}
Let \(0 < p < 1\,,\) \(0 < q \leq \infty \). Then\\

\noindent
(i) \((H^{p,q}, {\cal B}) = (H_0^{p,\infty }, {\cal B}) = H(\infty, \infty,1,1/p) =(H^p)^*\,,\)\\
(ii)\((H^{p,q}, {\cal B}_0) = (H_0^{p,\infty }, {\cal B}_0) = H(\infty, \infty,1,1/p) =(H^p)^*\,,q \not= \infty \,,\)\\
(iii)\((H^{p,\infty }, {\cal B}_0) = H_0(\infty, \infty,1,1/p)\,.\)\\ 
\end{Cor}

\noindent
We observe here that the fractional derivative operator \(D = D^1\) is a continuous isomorphism of 
the analytic Zygmund space \(\Lambda _\ast ^\infty (\mathbb{D})\) (resp. \(\lambda _\ast ^\infty (\mathbb{D})\)) onto \({\cal B}\) (resp. \({\cal B}_0\))
Thus, Corollary 8.2 may be viewed as an isomorphic version of the Zygmund space portion of Corollary 8.1.

In contrast, the Lipschitz spaces \(\Lambda _1^\infty (\mathbb{D})\) will, in general, determine different multiplier spaces 
\((H^{p,q},\Lambda _1^\infty (\mathbb{D}))\) for different values of \(q\). This is demonstrated in the next result which follows from Theorem 8.1 
and Theorem 4.2.


\begin{Cor}
Let \(0 < p < 1\,,\) \(0 < q < \infty\). Then\\

\noindent
(i) \((H^{p,q},\Lambda _1^\infty (\mathbb{D})) = (H^{\frac{p}{1 + p},q})^*\,,\)\\
(ii)\((H^{p,\infty},\Lambda _1^\infty (\mathbb{D})) = (H_0^{\frac{p}{1 + p},\infty })^*\,.\)\\
\end{Cor}

\noindent
For \(1 \leq s < \infty \,,\) we have the following analogs of Corollaries 8.1 and 8.3.


\begin{Cor}
Let \(0 < \alpha , p < 1 \leq s < \infty \,,\) \(0 < q \leq \infty\). Then\\

\noindent
(i) \((H^{p,q},\Lambda _\alpha ^s(\mathbb{D})) = (H_0^{p,\infty },\Lambda _\alpha ^s(\mathbb{D})) = H(s,\infty ,1, \frac{1}{p} + \alpha ) = (H^{\frac{p}{1 + \alpha p}},H^s)\,,\)\\
(ii) For \(q \not= \infty \,,\)\\ 
\[(H^{p,q},\lambda _\alpha ^s(\mathbb{D})) =  (H_0^{p,\infty },\lambda _\alpha ^s(\mathbb{D})) = H(s,\infty ,1, 1/p + \alpha )=(H^{\frac{p}{1 + \alpha p}},H^s)\,,\]\\
\noindent 
(iii) \((H^{p,\infty },\lambda _\alpha ^s(\mathbb{D})) = H_0(s,\infty ,1,\frac{1}{p} + \alpha )\,,\)\\
(iv) \((H^{p,q},\Lambda _\ast ^s(\mathbb{D})) = (H_0^{p,\infty },\Lambda _\ast ^s(\mathbb{D})) = H(s,\infty ,1, \frac{1}{p} + 1) = (H^{\frac{p}{1 + p}},H^s)\,,\)\\
(v) For \(q \not= \infty \,,\)\\ 
\[(H^{p,q},\lambda _\ast ^s(\mathbb{D})) =  (H_0^{p,\infty },\lambda _\ast ^s(\mathbb{D})) = H(s,\infty ,1, 1/p + 1)=(H^{\frac{p}{1 + p}},H^s)\,,\]\\
(vi) \((H^{p,\infty },\lambda _\ast ^s(\mathbb{D})) = H_0(s,\infty ,1,\frac{1}{p} + 1)\,.\)\\
\end{Cor} 


\begin{Cor}
Let \(0 < p < 1 \leq s < \infty \,,\) \(0 < q \leq \infty\). Then\\

\noindent
(i)\((H^{p,q},\Lambda _1^s(\mathbb{D})) = (H^{\frac{p}{1 + p},q},H^s)\,,\)\\
(ii)\((H_0^{p,\infty },\Lambda _1^s(\mathbb{D})) = (H_0^{\frac{p}{1 + p},\infty },H^s)\,.\)\\
\end{Cor}

\noindent
The space \(BMOA\) is the space of functions \(f \in H(\mathbb{D})\) having non-tangential limits \(m\)-a.e. on
\(\mathbb{T}\) for which the resulting boundary value function \(f\) is of bounded mean oscillation on \(\mathbb{T}\).
That is for which

\begin{equation}\label{8.1}
\sup_{I\subset\mathbb{T}}[m(I)^{-1}||(f - f_I)\chi _I||_1] < \infty  
\end{equation}

\noindent
where the supremum in \eqref{8.1} is taken over all subintervals \(I\subset\mathbb{T}\) and\\
\(f_I = m(I)^{-1}\int _I  f(z)\,dm(z)\). The space \(BMOA\) is not a Bergman-Sobolev space. However we do have the following embedding, see [38].


\begin{Lem}
\(H(\infty ,2,1,1)\hookrightarrow BMOA\).
\end{Lem}


\begin{Thm}
Let \(0 < p < 1\,,\) \(0 < q \leq 2\). Then
\[(H^{p,q},BMOA) = (H^p)^* = (H^{p,\infty },\cal B)\,.\]
\end{Thm}

\begin{proof}
Using the Duren-Romberg-Shields Theorem, Theorem 6.2, Lemma 8.1, Lemma 3.1(ii), and
Corollary 8.2(i) we find
\((H^p)^* = H(\infty,\infty,1,1/p) = (H^{p,q},H(\infty ,2,1,1) \subset (H^{p,q},BMOA)\subset (H^{p,q},{\cal B}) =
(H^p)^*\) and \((H^p)^* = (H^{p,\infty },\cal B)\) by Corollary 8.2(i). 
\end{proof}

We have not been able to find \((H^{p,q},BMOA)\) for \(2 < q \leq \infty \). However we can show 

\begin{equation}\label{8.2}
(H^{\frac{1}{2},\infty },BMOA) \not= (H^{\frac{1}{2},\infty },\cal B).
\end{equation}

\noindent
To show \eqref{8.2} let \(G\) be a function in \(\cal B\) which is not in \(BMOA\). The fractional
integral operator \(D = D_1\) is a continuous isomorphism of \(\cal B\) onto \(\Lambda _\ast ^\infty (\mathbb{D})\)
and hence \(G_{[1]} \in \Lambda _\ast ^\infty (\mathbb{D})\), [10]. But \(\Lambda _\ast ^\infty (\mathbb{D}) = (H^{\frac{1}{2},\infty },\cal B)\) by
Theorem 8.1 and Corollary 8.2(i) and so \(G_{[1]} \in (H^{\frac{1}{2},\infty },\cal B)\). Since the function \(f(z) = (1 - z)^{-2} \in H^{\frac{1}{2},\infty }\)
and since \(f\ast G_{[1]} = G\) then \(G_{[1]}\) fails to multiply \(H^{\frac{1}{2},\infty }\) into \(BMOA\). Thus \(G_{[1]}\) belongs to
\((H^{\frac{1}{2},\infty },\cal B)\) but not to  \((H^{\frac{1}{2},\infty },BMOA)\) which establishes \eqref{8.2}.

\end{document}